\documentclass[reqno]{amsart}
\usepackage[all]{xy}
\usepackage{amsmath}
\usepackage{amssymb,latexsym}
\usepackage{eucal}

\newtheorem{theorem}{Theorem}[section]
\newtheorem{proposition}[theorem]{Proposition}
\newtheorem{lemma}[theorem]{Lemma}

\newtheorem{definition}[theorem]{Definition}
\newtheorem{remark}[theorem]{Remark}
\newtheorem{example}[theorem]{Example}

\numberwithin{equation}{subsection}

\def\bee{\begin{equation}}
\def\lam{\lambda}
\def\lap{\lambda'}
\def\ga{\gamma}

\def\be{\beta}
\def\al{\alpha}
\def\de{\delta}
\def\te{\vartheta}
\def\cc{{\mathcal{C}}}
\def\ddc{\mathcal{D}}
\def\pc{{\mathcal{P}}}
\def\qc{{\mathcal{Q}}}

\def\gc{{\mathcal{G}}}
\def\vc{\mathcal{V}}

\def\prb{{p_{1}^{\ast}}}
\def\prc{{p_{2}^{\ast}}}

\def\la{\longrightarrow}

\def\R{{\mathbb{R}}}

\renewcommand{\include}{\input}
\begin{document}

\title{Notes on 1- and 2-gerbes}
\author{Lawrence Breen}

\address{UMR CNRS 7539 \newline
\indent Institut Galil\'ee \newline
\indent Universit\'e Paris 13 \newline
\indent  93430  Villetaneuse, France}
  \email{breen@math.univ-paris13.fr}

\maketitle

\begin{center}

\end{center}

\bigskip

The aim of these  notes is to  discuss in an   informal  manner 
 the construction and some
 properties of 1- and 2-gerbes. 
They are for the most part  based
 on the author's texts
\cite{festschrift}-\cite{2-gerbe}.
Our main goal is to describe the  construction  which associates to
 a gerbe or a 2-gerbe the corresponding 
 non-abelian cohomology class.

\bigskip

We begin by reviewing the well-known theory for principal
 bundles and show how to extend this
 to biprincipal bundles (a.k.a  bitorsors). After reviewing
 the definition  of stacks and gerbes, we construct
 the cohomology class associated to a gerbe. 
While the construction
 presented  is
 equivalent to that in \cite{2-gerbe}, it is clarified here
 by making  use of  
 diagram \eqref{verif:2coc}, a definite improvement over the corresponding
 diagram  \cite{2-gerbe}  (2.4.7),
 and  of \eqref{diagcoboun}. After a short discussion regarding the role of
 gerbes in algebraic topology, we pass from $1-$ to $2-$gerbes.
The construction of the associated cohomology classes follows the same
lines as for 1-gerbes, but with the additional degree of complication
entailed by passing from 1- to 2-categories, so that  it now involves
diagrams reminiscent of those in \cite{dgg}. Our emphasis will be on
 explaining how the fairly elaborate  equations which define cocycles
 and coboundaries may be reduced to terms which can be described
in  the tradititional formalism  
of non-abelian cohomology.

\bigskip

  Since the concepts discussed here are very
general, we have at times  not  made explicit the mathematical 
objects to which  they apply. For example, when we refer to 
``a space'' this might
mean a topological space, but also ``a scheme'' when one prefers to work
in an algebro-geometric context, or even  ``a sheaf''
and we place ourselves implicitly  in the category of such spaces,
schemes,
 or sheaves. Similarly, the standard notion of an $X$-group scheme $G$
 will correspond in a topological context to that of a bundle of
 groups on a space $X$.  By this we mean a total space $G$ above a
 space $X$ that is a group in the cartesian monoidal category of
 spaces over $X$. In particular, the fibers $G_x$ of $G$ at points $x
 \in X$ are topological groups, whose group laws vary continuously
 with $x$.

\bigskip

  Finally, in computing cocycles
 we will  consider  spaces $X$
endowed with a covering $\mathcal{U}:= (U_i)_{i\in I}$ by open sets, but the
discussion will remain valid 
when  the disjoint union $\coprod_{i \in I} \, U_i$ is replaced by an
arbitrary
 covering 
morphism $Y \la X$ of $X$  for a given  Grothendieck topology.
The emphasis in vocabulary will be  on
spaces rather than schemes, and we  have  avoided any
non-trivial result from algebraic geometry. In that sense, the  text is
implicitly directed towards topologists and category theorists rather
than algebraic  geometers, 
even though  we have not sought to make precise 
the category of spaces in which we work.

\bigskip

It is a pleasure  to thank Peter May,  Bob Oliver and Jim Stasheff  for their
 comments, and  William
Messing  for his
very careful reading  of a 
preliminary version of this text.

{\Large
\section{\bf Torsors and bitorsors}
\label{tbt}
}

\subsection
\noindent Let $G$ be a bundle of groups on a space  $X$. The following
definition of a principal space is standard, but note the occurence
of a structural bundle of groups, rather than simply  a constant one.
 We are in effect  giving ourselves a family of groups $G_x$, parametrized by
points
 $x \in X$,  acting
principally on the  corresponding  fibers $P_x$  of $P$.

\bigskip

\begin{definition}

 A left principal $G$-bundle (or
 left $G$-torsor) on a topological space  $X$ is a space $P \stackrel{\pi}{\la} X$ above $X$,
 together with
 a  left group action $G \times_X P \la P$ such that the induced morphism

\begin{equation}
\label{deftorsor}
\begin{array}{ccl}
 G  \times_X   P  & \simeq & P \times_X P \\
  (g,p) & \mapsto & (gp,\,p)
\end{array}
\end{equation}
is an isomorphism. We require in addition  that there exists 
a family of 
local sections $s_i: U_i \la P$, for some open cover $\mathcal{U}= 
(U_i)_{i \in I}$ of $X$.
The groupoid of left $G$-torsors on $X$ will be  denoted $\mathrm{Tors}(X,\,G)$.
\end{definition}

\noindent The choice of a family of local sections 
 $s_i:U_i \la P$, 
 determines
 a $G$-valued 1-cochain  $g_{ij}:U_{ij} \la G$, 
 defined above $U_{ij}:= U_i \cap U_j$
 by the equations
\begin{equation}
\label{1coch}
s_i = g_{ij} \, s_{j}           \quad \quad \forall\: i,j \in I \,.
\end{equation}
The $g_{ij}$  therefore satisfies the  1-cocycle
equation 
\begin{equation}
\label{1cocy}
g_{ik}= g_{ij}\, g_{jk}  
\end{equation}
above $U_{ijk}$. Two such  families of local sections $(s_i)_{i \in
  I}$  
and    $(s'_i)_{i \in I}$ on the same open cover $\mathcal{U}$  differ by  a 
  $G$-valued  0-cochain $(g_i)_{i \in I}$
 defined by
\begin{equation}
\label{newsec}   
      s'_i =g_i s_i  \qquad \qquad \forall i\in I
\end{equation}
 and for which the corresponding 1-cocycles $g_{ij}$ and   $g'_{ij}$
 are related to each other by  the coboundary relations
\begin{equation}
\label{1coboun}
g'_{ij} = g_i\, g_{ij} \, g_j^{-1}
\end{equation}
 This
 equation determines an equivalence relation  on the set of
1-cocycles $Z^1(\mathcal{U},\, G)$
\eqref{1cocy}, and the induced set of equivalence classes for this
equivalence relation is denoted  $H^1(\mathcal{U}, G)$. Passing to
the limit over open covers $\mathcal{U}$ of $X$ yields the
 \v{C}ech non-abelian cohomology set $\Check{H}^1(X,\, G)$, 
which  classifies
 isomorphism classes of $G$-torsors on $X$. This set
 is endowed with  a distinguished element,
 the class of the trivial left  $G$-torsor $T_G:= G \times X$.

\begin{definition}
Let $X$ be a space, and $G$ and $H$ a pair of bundles of groups on $X$. 
A  $(G,H)$-bitorsor  on $X$ is a space $P$ over $X$, together with 
fiber-preserving left and right actions of $G$ and $H$  on $P$, which commute
with each other and which define both a left $G$-torsor and a right
$H$-torsor structure  on $P$.  For any bundle of groups $G$,  a
$(G,G)$-bitorsor is
simply called  a $G$-bitorsor.
\end{definition}

\noindent A family of local 
sections $s_i$ of  a $(G,H)$-bitorsor  $P$ determines a local 
identification  of $P$ with both the trivial left $G$-torsor and the trivial
right $H$-torsor.  It therefore defines a family of local
isomorphisms $u_i: H_{U_i} \la  G_{U_i}$
between the restrictions  above  $U_i$
 of the bundles $H$ and $G$,
 which are explicitly given by the
rule
\begin{equation}
\label{defui0}
s_i h=  u_i(h)s_i
\end{equation}
for all $h \in H_{U_i}$.
 This however does not imply that
  the bundles of groups  $H$ and $G$  are globally isomorphic.

\bigskip

\noindent
\begin{example}
\rm  
 {\it i}) The trivial $G$-bitorsor on $X$: 
the right action of $G$ on the left $G$-torsor
$T_G$ is the trivial one, given  by
fibrewise right translation.  This bitorsor will also be   denoted $T_G$.

\medskip

\noindent \hspace{2.3cm}  {\it ii}) The group $P^{\text{ad}}:=
\text{Aut}_G(P)$ of $G$-equivariant fibre-preserving automorphisms of a
left $G$-torsor
$P$ acts on the right on $P$ by the rule
\[pu:= u^{-1}(p) \] 
so that any left $G$-torsor  $P$ is actually a
 $(G, P^{\text{ad}})$-bitorsor.  The group
  $P^{\text{ad}}$ is 
know as the gauge group of $P$. In particular, a left $G$-torsor $P$ 
 is a $(G,H)$-bitorsor if 
and only if the bundle of groups
$P^{\text{ad}}$  is isomorphic to $H$.

\medskip

\noindent \hspace{2.3cm}  {\it iii}) Let 
\begin{equation}
\label{sexs}
 1 \la G   \stackrel{i}{\la} H \stackrel{j}{\la} K \la 1 
\end{equation}
be a short exact sequence of bundles of groups on $X$.
 Then $H$ is a $G_K$-bitorsor on
$K$,  where
 the left and right actions above $K$ of the bundle of groups
$G_K:= G \times_X K$ 
   are  given by left and right
multiplication in $H$:
\[
\begin{array}{ccc}
(g,\,k) \ast h := f(g)\,h & \qquad & h \ast (g,\, k):= h\,f(g)
\end{array}\]
where $j(h)= k$.

\end{example}

\bigskip

\noindent

\subsection
\noindent Let  $P$ be a $(G, H)$-bitorsor and $Q$ be an
$(H,K)$-bitorsor
 on $X$.
 Let us define the  contracted product of $P$ and $Q$ as follows:

\begin{equation}
\label{def:P-twist}
 P \wedge^{H} Q : =   \frac{ P \times_X Q}{(ph, q) \sim   (p, hq)} 
\end{equation}
\noindent It is 
 a $(G,K)$-bitorsor on $X$ via the action of $G$ on $P$ and the
 action of  $K$ on $Q$. To any $(G,H)$-bitorsor $P$ on
$X$ is
associated the opposite $(H,G)$-bitorsor $P^{o}$, with same underlying
space as $P$, and for which the right action of $G$
( {\it resp.} left
action of $H$) 
 is induced by the given 
 left $G$-action ({\it resp.} right $H$-action)
 on $P$.
For a given bundle of groups $G$ on $X$, the category $\mathrm{Bitors}(X,\,G)$ 
 of $G$-bitorsors  on $X$ is a group-like  monoidal 
groupoid, in other words a monoidal category which is a groupoid, and 
in which every  object has
both a left and a right inverse. The tensor multiplication in 
 $\mathrm{Bitors}(X,\,G)$ 
 is
the contracted product of $G$-bitorsors,
 the unit object is the trivial bitorsor $T_G$, 
 and $P^{o}$ is an  inverse
 of the $G$-bitorsor  $P$. 
Group-like monoidal groupoids are also known as $gr$-categories.

\vspace{.5cm}

\subsection
\noindent {\bf Twisted objects:}

\vspace{.3cm}

Let $P$ be a left $G$-torsor on $X$, and  $E$ a space over  $X$ on which 
$G$ acts on the right. We say that  the space
$E^{P} :=E \wedge^G P$ over $X$,
 defined as in \eqref{def:P-twist}, is the $P$-twisted
form of $E$. The choice of a local section $p$ of $P$ above an open set $U$
determines an isomorphism $\phi_p: E^P_{|U}  \simeq E_{|U}$.
Conversely, if $E_1$ is a space over $X$ for which there exist a open cover
 $\mathcal{U}$ of $X$ above  which $E_1$ is 
locally isomorphic to $E$, then the space $\text{Isom}_X(E_1,\, E)$ is 
a left torsor on $X$ under the action of the bundle of groups 
$G:= \text{Aut}_X E$.

\begin{proposition}
These two constructions are inverse to each other.
\end{proposition}

\begin{example} {\rm  Let $G$ be a bundle of groups on $X$ and $H$ a bundle
of groups locally isomorphic to $G$
and let $P: = \mathrm{Isom}_X(H,\, G)$  be the left
$\mathrm{Aut}(G)$-torsor
  of fiber-preserving isomorphisms from $H$ to $G$. 
The map
\[
\begin{array}{ccl}
G \wedge^{\mathrm{Aut}(G)} P & \stackrel{\sim}{\la} & H\\
(g,\, u) & \mapsto & u^{-1}(g)
\end{array}
\]
identifies $H$ with the $P$-twisted form
of $G$, for the right action of $\mathrm{Aut}(G) $ on $G$ induced by
the standard left action.  Conversely, given a fixed bundle of groups
$G$ on $X$,  a $G$-torsor $P$ determines a bundle of
groups $H:= G \wedge^{\mathrm{Aut}(G)} P$ on $X$ locally isomorphic to $G$,
and $P$ is isomorphic to the left $\mathrm{Aut}(G)$-torsor
$\mathrm{Isom}(H,\,G)$.
}
\end{example}

\noindent The next  example is very well-known, but deserves to
be spelled out   in some
detail. 

\begin{example} 
{\rm
 A  rank
 $n$ vector bundle $\vc$ on $X$ is
locally isomorphic to the trivial bundle  $\mathbb{R}^n_X:=
X \times \mathbb{R}^n$, whose 
group of automorphisms is the trivial bundle of groups  
 \[GL(n,\mathbb{R})_X:=
GL(n,  \mathbb{R})
 \times X\]
on $X$. The  left principal  $GL(n,\,\mathbb{R})_X$-bundle associated to $\vc$ 
is its bundle of frames
$P_{\vc}
:= \text{Isom}(\vc,\,\mathbb{R}^n_X)$. 
The vector bundle
$\vc$ 
may be recovered from 
$P_{\vc}$ via the isomorphism 
\begin{equation}
\label{lin} \begin{array}{ccc}
 \mathbb{R}^n_X \wedge^{G L(n,{\mathbb{R}})_X}\,
 P_{\vc} &\stackrel{\sim}{\la} & \vc\\
( y,   p)  & \mapsto & p^{-1}(y)  
\end{array}  \end{equation}
\noindent in other words as the $P_{\vc}$-twist of the trivial vector bundle
$\mathbb{R}^n_X$ on $X$. Conversely, for any principal
$GL(n,\,\mathbb{R})_X$-bundle $P$ on $X$, the twisted object
 $\vc :=  \mathbb{R}^n_X \wedge^{G L(n,{\mathbb{R}})}\,
 P$ is known as the  rank $n$ vector bundle associated to
 $P$. Its frame bundle $P_{\vc}$ is canonically isomorphic to $P$.
}
\end{example}

\begin{remark}{\rm
In \eqref{lin}, the right action  on $\R^n_X$  of the  linear group
 $\mathrm{GL}(n,\,\R)_X$   is given by the rule
\[ \begin{array}{ccc} 
\R^n  \times \mathrm{GL}(n,\,\R) & \la & \R^n \\
(Y,\, A) & \mapsto & A^{-1} Y 
\end{array} 
\]
where an element of $\R^n$ is viewed as a column matrix $Y=
(\lambda_1, \ldots , \lambda_n)^T$. A local section $p$
of $P_{\vc}$ determines  a local basis $\mathcal{B} = \{p^{-1}(e_i)\}$
 of  $\vc$ and the arrow \eqref{lin} 
then identifies the column
  vector $Y$ with the 
element of $\vc$   with coordinates $(\lambda_i)$  in the chosen basis
$p$.
 The fact that the arrow
\eqref{lin}
 factors through the contracted product is a global version of
the familiar linear algebra rule which in  an $n$-dimensional  vector space $V$
 describes  the effect of  a change of basis matrix $A$ 
 on the  coordinates $Y$ of a given
vector $v \in V$. } 
\end{remark}

\bigskip

\subsection
\noindent{\bf The cocyclic description of a bitorsor
 (\cite{cg}, \cite{festschrift}):} \label{monoid}

\vspace{.3cm}

Consider a  $(G,H)$-bitorsor
$P$ on  $X$, with chosen  local sections  $s_i:U_i \la P$
 for some open cover $\mathcal{U}= (U_i)_{i \in I}$. Viewing $P$ as
a left $G$-torsor, we know by \eqref{1coch}   that these sections define
 a family of  $G$-valued  1-cochains  $g_{ij}$ 
 satisfying the 1-cocycle condition \eqref{1cocy}.
We have also seen that the  right  $H$-torsor
 structure on $P$ is then described by the family  of local
 isomorphisms  $u_i: H_{\,U_i} \la  G_{\,U_i}$ 
 defined by the equations \eqref{defui0}
for all $h \in H_{U_i}$.
It follows from  \eqref{1coch} and \eqref{defui0} that 
the transition law for the restrictions of these isomorphisms
above
  $U_{ij}$ is
\begin{equation}
\label{transu}
 u_i = i_{g_{ij}} \, u_j
\end{equation}
with $i$ the inner conjugation homomorphism 
\begin{equation}
\label{def:inconj}
\begin{array}{ccc}
 G & \stackrel{i}{\la} & \mathrm{Aut}(G) \\
g & \mapsto &i_g
\end{array} 
\end{equation}  
defined by 
\begin{equation}
\label{inner0}
i_g(\ga)= g \ga g^{-1}\,.
\end{equation}
The pairs $(g_{ij},\, u_i)$ therefore satisfy the cocycle conditions
 \begin{equation}
\label{defbitorcoc}
\begin{cases}
 g_{ik}  =\,\, \, g_{ij}\, g_{jk}\\
u_i = \,\,\, i_{g_{ij}}\, u_j
\end{cases}
\end{equation}
 A  second family of local sections $s'_i$ of $P$
 determines a corresponding cocycle pair $(u'_i,\,g'_{ij})$,
 These new cocycles differ from the previous ones
by the coboundary relations 
\begin{equation}
\label{bitcoboun}
\begin{cases}
g'_{ij} = g_i\, g_{ij} \, g_j^{-1} \\
u'_i  = i_{g_i} \,u_i 
\end{cases}
\end{equation}
where  the 0-cochains $g_i$ are  defined by \eqref{newsec}.
  Isomorphism classes of $(G,H)$-bitorsors on $X$  with given
local trivialization on an open covering $\mathcal{U}$ are
classified  by the quotient   of the set of  
 cocycles $(u_i,\,g_{ij})$    \eqref{defbitorcoc}
by  the equivalence relation  \eqref{bitcoboun}. 
Note that when $P$ is a $G$-bitorsor,  the terms of the  second
equation
in both \eqref{defbitorcoc} and \eqref{bitcoboun}
live in the group $\mathrm{Aut}(G)$.   
In that case, the set of cocycle classes is  the non-abelian
hypercohomology set
$H^0(\mathcal{U},\, G \la \mathrm{Aut}(G))$, with values in the length
one complex
of groups \eqref{def:inconj} where $G$ is placed in degree $-1$.
 Passing to the limit over open
covers, we obtain the \v{C}ech cohomology set
  $\Check{H}^0(X, \, G \la \mathrm{Aut}(G))$ which classifies
  isomorphism classes of $G$-bitorsors on $X$.

\bigskip

 Let us see how the  monoidal structure on the category
 of $G$-bitorsors is reflected at the cocyclic level.  Let
 $P$ and $Q$ be a pair of $G$-bitorsors on $X$, with chosen local
 sections $p_i$ and $q_i$. These determine  corresponding cocycle
 pairs  $(g_{ij},\, u_i)$ and $(\ga_{ij},\,v_i)$ satisfying the
 corresponding equations \eqref{defbitorcoc}. It is readily verified
 that
the corresponding
 cocycle pair for the $G$-bitorsor $P \wedge^G Q$, locally trivialized
 by the family of local sections $p_i \wedge q_i$, is the pair
\begin{equation}
\label{prod}
(g_{ij}\, u_i(\ga_{ij}),\, u_i\, v_i)
\end{equation}
so that the group law for cocycle pairs is simply
 the semi-direct product multiplication in the group  $G \rtimes
\mathrm{Aut}(G)$, for the standard left action of $\mathrm{Aut}(G)$
on $G$. The multiplication rule for  cocycle pairs
\[(g_{ij},\, u_i) \ast (\ga_{ij},\,v_i) =
 (g_{ij}\, u_i(\ga_{ij}),\, u_i\, v_i) \] passes to the set of
 equivalence classes,   and  therefore 
 determines a group structure on the set \linebreak
  $\Check{H}^0(X, \, G \la
 \mathrm{Aut}(G))$, which reflects the contracted product of bitorsors.

\begin{remark}
\label{inner}
{\rm Let us choose once more  a family of local sections $s_i$ of a
  $(G,H)$-bitorsor $P$. The local isomorphisms $u_i$  provide an
 identification
of   the restrictions $H_{U_i}$ of $H$ with the restrictions  $G_{U_i}$ of
  $G$. Under these identifications, the significance of  equations
  \eqref{transu} is the following. By \eqref{transu}, 
we may  think of an element of
  $H$ as given by a family  of local elements $\ga_i \in
  G_i$,  glued to each other
above  the open sets $U_{ij}$ according to the rule
\[ \ga_i =i_{g_{ij}}\, \ga_j \,. \]
For this reason, a bundle of  groups $H$ which stands
 in such a relation to a
given group $G$ may be  called an {\it inner form} of 
$G$. This is the terminology used in the context of Galois cohomology, 
 {\it
  i.e}
when $X$  is a scheme $\mathrm{Spec}(k)$ endowed with the
\'etale topology defined by  the covering morphism
 $\mathrm{Spec}(k') \,\la\,\mathrm{Spec}(k)$
 associated
 to a Galois field extension  $k'/k$ (\cite{cg} III \S 1). 
}
\end{remark}

\bigskip

\subsection \noindent  The previous  discussion remains  valid in a wider
context, in which  the inner conjugation homomorphism $i$
 is replaced by an
arbitrary homomorphism of groups $\delta: G \la \Pi$. The cocycle and
coboundary conditions \eqref{defbitorcoc} and \eqref{bitcoboun}
are now respectively  replaced by the rules
\begin{equation}
\label{defbitorcoc1}
\left\{\begin{array}{cl}  g_{ik} & =\,\, \, g_{ij}\, g_{jk}\\
\pi_i &= \,\,\, \de(g_{ij})\, \pi_j
\end{array}
\right.
\end{equation} 
and by
\begin{equation}
\label{bitcoboun1}
\begin{cases}
g'_{ij} = g_i\, g_{ij} \, g_j^{-1} \\
\pi'_i  = \de (g_i) \,\pi_i 
\end{cases}
\end{equation}
and the induced 
  \v{C}ech hypercohomology set  with values in the complex of groups $G
\la \Pi$ is denoted
$\Check{H}^0(\mathcal{U}, G \la \Pi)$.
In order to extend to $\Check{H}^0(\mathcal{U}, G \la \Pi)$
 the multiplication \eqref{prod}, we require
additional structure:

\begin{definition}
A (left) crossed module is a group homomorphism $\de:  G \la \Pi$,
together with a left group action
\[ \begin{array}{ccc}
\Pi \times G &\la &G\\ 
(\pi,\,g) & \mapsto & {}^{\pi\!}g
\end{array}\]
of $\Pi$ on the group $G$, and  such that the equations
\begin{equation}
\begin{cases}
\de ( {}^{\pi\,}\!g ) =  {}^{\pi\,}\!\de(g)\\
 {}^{\de(\ga)\,}\!g  =  {}^{\ga\,}\!g
\end{cases}
\end{equation} 
are satisfied, with $G$ ({\it resp.} $\Pi$)  acting on itself by the 
 conjugation rule \eqref{inner0}.
\end{definition} 
Crossed modules form a category, with a homomorphism of crossed modules
\[(G \stackrel{\de}{\la} \pi) \la ((K \stackrel{\de'}{\la} \Gamma)\]
defined by a pair of homomorphisms $(u,v)$ such that the diagram of
groups
\begin{equation}
\label{crmodhom}
\xymatrix@C=30pt{
G \ar[r]^u \ar[d]_{\de} & K \ar[d]^{\de'} \\
\Pi \ar[r]_{v} & \Gamma
}\end{equation}
commutes, and such that
 $u({}^{\pi\,}\!g) = {}^{v(\pi)\,}\! u(g)$ (in other words such that 
$u$ is $v$-equivariant). 

\bigskip

A left crossed module $G \stackrel{\de}{\la} \Pi$ defines 
 a group-like monoidal
category $\cc$ with a strict multiplication on objects, by setting
\begin{equation}
\label{defmon}
\mathrm{ob}\, \cc := \Pi \qquad \mathrm{ar}\,\cc:= G \times \Pi
\end{equation}
The source and target of an arrow $(g,\, \pi)$ are  as follows:
\[\xymatrix@C=35pt{
\pi \ar[r]^(.45){(g,\pi)} & \de(g) \pi 
}\]
and the composite of two  composable arrows 
\begin{equation}
\label{comprule}
\xymatrix@C=55pt{
\pi \ar[r]^{(g,\pi)} & \de (g)\pi \ar[r]^(.45){(g',\, \de (g) \pi)} &
\de (g'g), \pi
}
\end{equation}
is the arrow $(g'g,\, \pi)$. 
The monoidal structure on this groupoid
is given on the objects by the group multiplication in $\Pi$, and on the
set $G \times \Pi$ of arrows by the semi-direct product group multiplication
\begin{equation}
\label{semi1}
(g,\, \pi) \ast (g'\, \pi') := (g\, \,{}^{\pi\,}\!g',\, \pi \, \pi')
\end{equation}
for the given left action of $\Pi$ on $G$.
In particular the identity element of the group $\Pi$ is the unit
object $I$ of this monoidal groupoid.

\bigskip

 Conversely, to a monoidal
category $\mathcal{M}$ with strict multiplication on objects is
associated a crossed module $G \stackrel{\de}{\la} \Pi$, where $\Pi: =
\mathrm{ob}\, \mathcal{M}$ and
 $G$ is the set $\mathrm{Ar}_I \mathcal{M}$ of arrows of 
$\mathcal{M}$ sourced at the identity object, with  $\de$
 the restriction to $G$ of the target map. The group law on $G$ is the
 restriction  to this set of the multiplication of
 arrows in  the monoidal category $\mathcal{M}$. The action of
 an object $\pi \in \Pi$ on  an arrow $g: I \la \de(g)$ in $G$
has the following categorical interpretation: the composite arrow
\[\xymatrix@C=40pt{I \ar[r]^(.4){\sim} &
 \pi\, I \,\pi^{-1} \ar[r]^{\pi\, g\,  \pi^{-1}}&
  \pi\, \de (g) \,\pi^{-1}} \]
 corresponds to  the element  ${}^{\pi\,}\!g$ in $G$.
Finally,
given a pair elements
 $g,g' \in  \mathrm{Ar}_I \mathcal{M}$, it
 follows from the composition rule   \eqref{comprule} for a  pair of arrows
 that the composite arrow
\[ \xymatrix@C=55pt{
I \ar[r]^{(g,I)} & \de(g) \ar[r]^{(g',\de(g))}& \de(g'g)
}\]
 (constructed by taking advantage of the monoidal structure on the
 category $\mathcal{M}$ in order to transform the arrow $g'$ into an arrow
$(g',\,\de(g))$  composable  with $g$) is simply given by the element $g'g$ of the group 
$ \Pi = \mathrm{Ar}_I \mathcal{M}$.

\bigskip

A stronger concept  than that of a homomorphism of crossed module is
what could be termed  a ``crossed module of crossed modules''. This is the
categorification of   crossed modules and corresponds,
when one extends the previous dictionary between strict monoidal categories
and crossed modules, to strict monoidal bicategories. The most
efficient description of such a concept  is  the notion of a crossed square, due
to J.-L Loday. This consists of a homomorphism of crossed
modules $\eqref{crmodhom}$, together with a map 
\begin{equation}
\label{crhommap}
\begin{array}{ccc}
K \times \Pi & \la& G\\
(k,\,\pi) & \mapsto & \{k,\, \pi\} 
\end{array}  
\end{equation} 
satisfying certain conditions for which we refer to 
 \cite{loday} definition 5.1.
 
\begin{remark}
\label{rem:crmod}
 {\rm {\it i)} 
The definition \eqref{crmodhom}
 of a homomorphism of crossed modules   is
quite restrictive, and it is often preferable to relax it so that
 it defines a not necessarily strict monoidal functor between the
 associated (strict) monoidal groupoids. The  definition of a weak
 homomorphism of crossed modules has been spelled out by
 B. Noohi  in
 \cite{noohi}  (definition 8.4), to which we also refer for a
 discussion of related issues.

\medskip 

\hspace{2cm}
{\it ii)}
 All these definitions obviously extend
 from groups to bundles of groups on
$X$.
\medskip

\hspace{2cm}
{\it iii)}
 The composition  law \eqref{semi1} determines a multiplication  
\[( g_{ij},\, \pi_i)  \ast (g'_{ij}, \pi'_i) :=
 (g_{ij}\,{}^{\pi_i\,}\!g'_{ij},\, \pi_i\, \pi_{i}') \]
 on $(G \la \Pi)$-valued cocycle pairs,
which generalizes \eqref{prod},  is compatible with the coboundary
relations,   and induces a group structure  on the  set
 $\Check{H}^0(\mathcal{U},\, G \la \Pi)$ 
of degree zero cohomology classes with
values in the crossed module $G \la \Pi$ on $X$.
\hspace{2cm}
}
\end{remark}

\subsection \noindent The following proposition 
 is known as the Morita theorem, by 
analogy with the corresponding characterization
in terms  of bimodules of  
equivalences between certain categories of modules.

\begin{proposition} \label{morita}  {\rm (Giraud \cite{gir})}
{\it i}) A $(G,H)$-bitorsor $Q$  on $X$ determines an equivalence 
\[ \begin{array}{ccc}
\mathrm{Tors}(H) & \stackrel{\Phi_Q}{\la}  & \mathrm{Tors}(G)\\
 M & \mapsto & Q   \wedge^{H} M\end{array}
\]
between the corresponding categories of left torsors on $X$. In
addition, if $P$ is an $(H,K)$-bitorsor on $X$, then there is a natural
equivalence
\[\Phi_{\,Q \,\wedge^{H}\,P} \simeq \Phi_Q \circ \Phi_P \]
 between functors from $\mathrm{Tors}(K)$ to $\mathrm{Tors}(G)$.
 In particular, the equivalence $\Phi_{Q^o}$ in an inverse of $\Phi_Q$.

\medskip

\hspace{4cm}{\it ii})  Any such equivalence $\Phi$ between two
 categories of torsors  is
equivalent to one associated in this manner to an $(H,G)$-bitorsor.
\end{proposition}

\bigskip

\noindent{\bf Proof of} $\boldsymbol{i}\boldsymbol{i}):$ To a given 
 equivalence
 $\Phi$ is associated the left $G$-torsor 
$Q:= \Phi(T_H)$. By functoriality of $\Phi$,
$H \simeq \text{Aut}_H(T_H) \stackrel{\Phi}{\simeq}  \text{Aut}_G(Q)$, so that  a section  of
$H$ acts   on the right on $Q$.

{\Large
\section{\bf (1)-stacks}
}
\subsection
\noindent The concept of a stack is the categorical analogue of a sheaf.
Let us  start by defining the analog of a presheaf.

\begin{definition}
{\it i)}:
  A  category fibered in groupoids above a space $X$ 
consists in a family of groupoids $\cc_U$, for each open set $U$ in $X$,
 together with an inverse image functor
\begin{equation}
\label{definv}
 f^{\ast}: \cc_U \la \cc_{U_1} 
\end{equation}

\noindent associated to every inclusion of open sets  $f: U_1 \subset U$
 (which is the identity
 whenever $f = 1_U$), and  natural equivalences

\begin{equation}
\label{defphi}
 \phi_{f,g}: (fg)^{\ast} \Longrightarrow g^{\ast}\, f^{\ast} 
\end{equation}
for every pair of composable inclusions
\begin{equation}
\label{compos}
 U_2 \stackrel{g}{\hookrightarrow} U_1  \stackrel{f}{\hookrightarrow}
 U \:.
\end{equation}

\noindent For each triple of composable inclusions 

\[ U_3 \stackrel{h}{\hookrightarrow}   U_2 \stackrel{g}{ \hookrightarrow}
 U_1  \stackrel{f}{ \hookrightarrow} U\,. \]

\noindent  we also
require  that the composite natural transformations 
\[ \psi_{f,g,h}: (fgh)^{\ast} \Longrightarrow h^{\ast}\,(fg)^{\ast}
 \Longrightarrow h^{\ast} \, (g^{\ast} f^{\ast}) \]
\noindent and
\[\chi_{f,g,h}:  (fgh)^{\ast} \Longrightarrow (gh)^{\ast}\, f^{\ast}
 \Longrightarrow (h^{\ast} g^{\ast}) \, f^{\ast}. \]
\noindent coincide.

\bigskip

\hspace{2cm} {\it ii)} A cartesian functor $F: \cc \la \ddc$
is a family of functors  $F_U: \cc_U \la \ddc_U$ for all open sets $U
\subset X$, together with natural transformations
\begin{equation}
\label{cart}
\xymatrix@R=15pt@C=30pt{\cc_U \ar[d]_{F_U} \ar[r]
 & \cc_{U_1} \ar[d]^{F_{U_1}}
  \ar@{}[d]_(.3){\,}="1"  \\
\ddc_U \ar[r] \ar@{}[r]_(.6){\,}="2"  & \ddc_{U_1}
 \ar@{}"2";"1"^(.2){\,}="3"
 \ar@{}"2";"1"^(.8){\,}="4"
\ar@{=>}"3";"4"
}\end{equation}
\noindent for all inclusion $f: U_1 \subset U$ compatible via the 
natural equivalences \eqref{defphi} for a pair of composable
inclusions \eqref{compos}

\bigskip

\hspace{2cm} {\it iii)} A natural transformation
 $\Psi:  F \Longrightarrow G$  between a pair of  cartesian
 functors consists of a family of natural transformations 
$\Psi_U: F_U \Longrightarrow G_U $ compatible via the 2-arrows
 \eqref{cart} under the inverse images functors \eqref{definv}.
\end{definition}

The following is the analogue for fibered groupoids of the notion of a 
sheaf of sets, formulated here 
in a preliminary style:

\begin{definition}
  A stack in groupoids  above a space $X$ is a fibered category
 in groupoids above $X$ such that
 \begin{itemize}
 \item (``Arrows glue'') For every pair of objects $x, y 
\in \cc_U$, the presheaf    $ \mathrm{Ar}_{\cc_U}(x,\,y)$ 
is a  sheaf on $U$.

\item (``Objects glue'') Descent is effective  for objects in  $\cc$. 
 \end{itemize}
\end{definition}

The gluing condition on arrows is not quite correct as stated. In
order to be more precise, let us first observe that if $x$ is any
object in $\cc_U$, and $(U_\alpha)_{\alpha \in I}$ an open cover of $U$,
then $x$   determines a family of inverse images $x_\alpha$  in
$\cc_{U_\alpha}$ which we will refer informally to as the restrictions
of $x$ above $U_\alpha$, and sometimes denote by  $x_{|U_\alpha}$.
These are endowed with isomorphisms
\begin{equation}
\label{descent1}
\xymatrix@C=40pt{
{x_\beta}_{|U_{\alpha \beta}} \ar[r]^{\phi_{\alpha \beta}} &
 {x_\alpha}_{|U_{\alpha \beta}}
  }
\end{equation}
in $\cc_{U_{\al \be}}$ satisfying the cocycle equation
\begin{equation}
\label{descent2}
\phi_{\al \be} \, \,\phi_{\be \ga} = \phi_{\al \ga}
\end{equation}
when restricted to $\cc_{U_{\al \be \ga}}$.
An arrow $f:x \la y$ in $\cc_U$ determines   arrows $f_\al:
x_\al \la y_\al$ in  each of the categories $\cc_{U_\al}$ such 
 that the following diagram  in $\cc_{U_{\al \be}}$ commutes
\begin{equation}
\label{descent3}
\xymatrix@C=40pt@R=30pt{
{x_\be}_{|U_{\al \be}} \ar[r]^{{f_\be}_{|U_{\al \be}}}
\ar[d]_{\phi_{\al \be}} & {y_\be}_{|U_{\al \be}}
\ar[d]^{\psi_{\al \be}}\\
{x_\al}_{|U_{\al \be}} \ar[r]^{{f_\al}_{|U_{\al \be}}} & {y_\al}_{|U_{\al \be}}
}
\end{equation}
The full gluing condition on arrows is the requirement that
conversely,
 for any family of arrows $f_\al: x_\al \la y_\al$ in $\cc_{U_\al}$
for which the compatibility condition \eqref{descent3} is satisfied,
there exists a unique arrow $f: x \la y$ in $\cc_U$ whose restriction
above each open set $U_\al$ is the corresponding 
 $f_\al$.  In particular, if we make the very
non-categorical additional  assumption that the $\phi_{\alpha \be}$ are all identity
arrows, then  this gluing condition on arrows  simply  asserts that
the presheaf of arrows from $x$ to $y$ is a sheaf on $U$.  A fibered
category for which the gluing property on arrows is satisfied is
called a prestack.

\bigskip

Let us now pass to  the gluing condition on objects. The term  descent
comes from algebraic geometry, where for a given  family of objects
$(x_\al) \in
\cc_{U_\al}$, 
a family of isomorphisms $\phi_{\al \be }$ \eqref{descent1} satisfying
the equation \eqref{descent2} is called descent data for the family of
objects $(x_\al)_{\al \in I}$. The descent is said to be effective
whenever any  such descent data determines an object  $x \in \cc_U$,
together with a family of arrows $x_{|{U_\al}} \la x_\al $ in
$\cc_{U_{\al}}$ compatible with the given descent data on the given
objects $x_\al$, and the canonical descent data on the restrictions
$x_{|U_\al}$ of $x$.
A   sheafification process, analogous to the one
which transforms a presheaf into a sheaf,  associates a stack to a
given prestack. For a more detailed introduction to the theory of
stacks in an algebro-geometric setting, see \cite{vistoli}.

{\Large
\section{\bf 1-gerbes}
}

\subsection
\noindent We begin with the  global description of the 2-category of  gerbes, due to
Giraud \cite{gir}. For another early discussion of gerbes, see 
\cite{duskin}.

\bigskip

\begin{definition}
i) A (1)-gerbe on a space $X$ is a stack in groupoids $\gc$ on $X$ which
 is locally {\bf  non-empty} and locally {\bf connected}.

\medskip
   \hspace{2cm}   ii) A morphism of gerbes  is a cartesian functor
 between the underlying stacks.

\medskip
  \hspace{2cm}   iii) A  natural 
 transformation   $\Phi:u \Longrightarrow v$  between a pair of  such
 morphisms
of gerbes $u,v: \mathcal{P} \la \mathcal{Q}$ is  a natural transformation
 between the corresponding pair of  cartesian functors.
\end{definition}

\begin{example}
{\rm  Let $G$ be a bundle of groups on $X$. The stack
$\cc := \text{Tors}(G)$  of left $G$-torsors on $X$ is a gerbe on $X$:
first of all, it is non-empty, since the category $\cc_U$ always has
at least one object, the trivial torsor $T_{G_U}$. In addition, every
$G$-torsor on $U$ is locally isomorphic to the trivial one, so the
objects in  the category $\cc_U$ are  locally connected. 
}
\end{example}

A gerbe $\pc$ on $X$  is said to be  {\it neutral} (or trivial) 
when the fiber category $\pc_X$ is non empty.
In particular, a gerbe $\text{Tors}(G)$ is neutral with distinguished object
the trivial $G$-torsor $T_G$ on $X$.
Conversely,
the choice of a global  object $x \in \pc_X$ in a neutral gerbe $\pc$ determines
an equivalence of gerbes
\begin{equation}
\label{globject}
\begin{array}{ccc}
\pc & \stackrel{\sim}{\la} & \text{Tors}(G)  \\
 y & \mapsto & \text{Isom}_{\pc}(y,\, x)
\end{array}
\end{equation}
 on $X$, where $G:= \text{Aut}_{\pc}(x)$, acting on
 $\mathrm{Isom}_{\pc}(x,\,y)$ by composition of arrows.

\bigskip

Let  $\pc$ be a gerbe on $X$ and  $\mathcal{U} = (U_i)_{i\in I}$ be an
 open cover of $X$.
 \noindent  We now {\bf choose} objects  $x_i \in \text{ob}\  \pc_{U_i}$
for each $i \in I$.
  These objects
 determine corresponding  bundles of groups $G_i:= \text{Aut}_{\pc_{U_i}}(x_i)$
 above $U_i$.
\noindent  When in addition there exists  a bundle of groups $G$ above $X$,
 together with $U_i$-isomorphisms $G_{|\,U_i} \simeq G_{i}$, for all
 $i \in I$,
 we say that $\pc$ is a
 $G$-gerbe on $X$.

{\Large
\section{\bf Semi-local description of a gerbe}
}
\label{semi}
\subsection
\noindent Let $\pc$ be  a $G$-gerbe   on $X$, and let us choose  a family of
 local objects $x_i \in \pc_{U_i}$. These determine as in \eqref{globject}
equivalences 
\[\Phi_i: \pc_{U_i} \la  \text{Tors}(G)_{|\, U_i}\]
above $U_i$.  Chosing quasi-inverses for the $\Phi_i$ we get an induced family of equivalences
\[\Phi_{ij} := \Phi_{i\, | U_{ij}} \circ  \Phi^{-1}_{j\, | U_{ij}}:
 \text{Tors}(G)_{U_{ij}}   \la \pc_{|U_{ij}} \la  \text{Tors}(G)_{U_{ij}} \]
\noindent above $U_{ij}$, which corresponds  by proposition \ref{morita}
to  a family of $G$-bitorsors $P_{ij}$ above $U_{ij}$. 
\noindent By construction of the $\Phi_{ij}$, there are also natural transformations
\[\Psi_{ijk}: \Phi_{ij} \, \Phi_{jk} \Longrightarrow \Phi_{ik}\]
 above $U_{ijk}$, satisfying a coherence condition on $U_{ijkl}$.
These define  isomorphisms of $G$-bitorsors
\begin{equation}
\label{tors:bitors}
\psi_{ijk}: P_{ij} \wedge^G P_{jk} \la  P_{ik}
\end{equation}
\noindent above $U_{ijk}$ for which 
 this coherence
condition is  described by the commutativity of the
 diagram of bitorsors
\begin{equation}
\label{bitorcocycle} \xymatrix@R=10pt@C=30pt{
P_{ij} \wedge P_{jk} \wedge P_{kl}\ar[dd]_{P_{ij}\wedge\psi_{ijk}}
\ar[rr]^(.55){\psi_{ijk}\wedge P_{kl}} && P_{ik} \wedge P_{kl}  \ar[dd]^{\psi_{ikl}} \\
&& \\
P_{ij} \wedge P_{jl} \ar[rr]_{\psi_{ijl}} && P_{il}
}\end{equation}
 above $U_{ijkl}$

\vspace{.5cm}

\subsection
\noindent {\bf Additional comments}:

\vspace{.3cm}

 {\it i}) 
The isomorphism \eqref{tors:bitors}, satisfying the coherence condition
 \eqref{bitorcocycle},
 may be viewed
as a 1-cocycle condition
 on $X$ with values in the monoidal stack of $G$-bitorsors
on $X$. We say that a family of such  bitorsors
 $P_{ij}$ constitutes  
   a bitorsor cocycle
on $X$.

\bigskip

 {\it ii}) In the case of { \it abelian} $G$-gerbes\footnote{which are
 not simply $G$-gerbes for which the structure group $G$ is abelian !}  
(\cite{2-gerbe} definition 2.9), the monoidal stack of bitorsors on $U_{ij}$ 
may be replaced by the
 symmetric monoidal stack of $G$-torsors on $U_{ij}$.
  In particular, for the multiplicative group  $G = GL(1)$, the
$GL(1)$-torsors 
 $P_{ij}$ correspond to  line bundles $L_{ij}$. This   the point of view regarding abelian $GL(1)$
gerbes set forth
 by N. Hitchin  in \cite{hitch}.

\bigskip

{\it iii}) The semi-local construction extends from
 $G$-gerbes to general gerbes. In that case a local group $G_i:=
 \mathrm{Aut}_{\pc}(x_i)$ above $U_i$ 
is associated to each of the chosen objects $x_i$. The previous
discussion remains valid, with the proviso that the 
  $P_{ij}$ are now 
$(G_j,\, G_i)$-bitorsors rather than simply $G$-bitorsors,
 and  the $\psi_{ijk}$  \eqref{tors:bitors}
are  isomorphisms of $(G_k,\,G_i)$-bitorsors.

\bigskip

{\it iv})  If we replace the 
chosen trivializing 
open cover $\mathcal{U}$ of $X$ by a single covering morphism $Y \la X$ in some 
Grothendieck topology, the theory remains unchanged, but takes on a
somewhat different flavor. An object $x \in \pc_Y$ determines a 
bundle of groups   $G:= \text{Aut}_{\pc_Y}(x)$ over $Y$,
 together with a $(\prc G,\prb G)$-bitorsor
$P$ above $Y\times_X Y$ satisfying the coherence condition analogous
to \eqref{bitorcocycle} above $Y \times_X Y\times_X Y$.
A bitorsor $P$ on $Y$ satisfying this coherence condition has been called
a cocycle bitorsor by K.-H. Ulbrich \cite{ulbrich}, and a bundle
 gerbe by M.K. Murray  \cite{murray}. It corresponds to a bouquet
 in  Duskin's theory (see \cite{ulbrich1}).
 It is equivalent\footnote{ 
For a more 
detailed discussion of this when the covering  morphism $Y \la X$  is
the morphism of schemes associated as in remark \ref{inner} to a Galois
field extension $k'/k$, see \cite{lb:tc} \S 5.} to give oneself such a
bundle gerbe $P$, or 
 to consider  a gerbe $\pc$ on $X$,
together with a trivialization of  its pullback   to $Y$, since to  a
trivializing object $x \in \mathrm{ob}(\pc_Y)$ we may associate the
$G$-bitorsor $P:= \mathrm{Isom} (p_2^\ast x,\, p_1^\ast x)$ above $Y
\times_X Y$.

{\Large
\section{\bf Cocycles and coboundaries for  gerbes}
}
\label{gerbecoc}
\subsection
 \noindent Let us keep the notations of  section \ref{semi}.  In
 addition to choosing local objects $x_i \in \pc_{U_i}$ in a gerbe
 $\pc$ on $X$,  we now 
{\bf  choose}     arrows 
\begin{equation}
\label{def:lij}
\xymatrix@=15pt{
x_j \ar[rr]^{\phi_{ij}} && x_i 
}\end{equation}
in  $\pc_{U_{ij}}$\footnote{\label{paracomp}
Actually, this is a simplification, since the gerbe axioms
only allow us to choose such an arrow locally, above each element
 $U_{ij}^{\alpha}$ of an open cover of $U_{ij}$. Such families of open
 sets $(U_i,\, U_{ij}^{\alpha})$, and so on,
 form what is known as a hypercover of $X$.
For simplicity, we assume from now on 
 that our topological
 space $X$ is 
paracompact. In that case,  we
may  carry
out the entire discussion without hypercovers.}.
Since   $G_i:=  \mathrm{Aut}_{\pc}(x_i)$, a chosen  arrow $\phi_{ij}$
 induces by conjugation a  homomorphism
 of group bundles
\begin{equation}
\label{def:conjgr}
\xymatrix@C=10pt@R=3pt{ G_{j\,|U_{ij}} \ar[rr]^{\lambda_{ij}}
 && G_{i\,|U_{ij}}
\\
\ga \ar@{|->}[rr] && \phi_{ij} \,  \ga \,\phi_{ij}^{-1}
}\end{equation}
\noindent above the open sets  $U_{ij}$. To state this slightly
differently, such a homomorphism $\lambda_{ij}$ is 
characterized by the commutativity of the diagrams
\begin{equation}
        \label{def-lamij}
        \xymatrix@R10pt@C=15pt{
   x_j \ar[rr]^{\ga} \ar[dd]_{\phi_{ij}}
    && x_j \ar[dd]^{\phi_{ij}}\\&&\\x_i 
    \ar[rr]_{\lambda_{ij}(\ga)}&& x_i
    }\end{equation} 
\noindent for every  $\ga \in G_{j\, |U_{ij}}$.
The choice of objects $x_i$ and arrows $\phi_{ij}$
 in $\pc$  determines, in addition to the morphisms $\lambda_{ij} $
 \eqref{def:conjgr},
 a family of
  elements  
 $g_{ijk} \in G_{i\,\mid U_{ijk}}$ for all $(i,j,k)$, defined 
by the commutativity of the  diagrams
\begin{equation}
        \label{coc0}
        \xymatrix@R=10pt@C=15pt{
   x_k \ar[rr]^{\phi_{jk}} \ar[dd]_{\phi_{ik}}
    && x_j \ar[dd]^{\phi_{ij}}\\&&\\x_i 
    \ar[rr]_{g_{ijk}}&& x_i
    }\end{equation} 
\noindent above $U_{ijk}$.
These in turn induce by conjugation the following
  commutative diagrams of bundles of  groups 
\begin{equation}
        \label{coc1}
        \xymatrix@R=10pt@C=15pt{
   G_k \ar[rr]^{\lambda_{jk}} \ar[dd]_{\lambda_{ik}}
    && G_j \ar[dd]^{\lambda_{ij}}\\&&\\G_i 
    \ar[rr]_{i_{g_{ijk}}}&& G_i
    }\end{equation} 
\noindent above $U_{ijk}$. 
\noindent The commutativity of diagram \eqref{coc1} may be
 stated 
  algebraically as
the cocycle equation
\begin{equation}
  \label{eq:coc1}
  \lambda_{ij} \, \lambda_{jk} = i_{g_{ijk}} \, \lambda_{ik}
\end{equation}
with $i$ the inner conjugation arrow \eqref{def:inconj}.
 The following equation is the second cocycle equation satisfied by the
pair $(\lambda_{ij},\, g_{ijk})$. While the proof of lemma
\ref{lem:coc2} given here is essentially the same as the one in
\cite{2-gerbe}, the present cubical diagram \eqref{verif:2coc}  is much  more
intelligible than diagram (2.4.7) of \cite{2-gerbe}.
\begin{lemma}
\label{lem:coc2}
The  elements $g_{ijk}$ satisfy the $\lambda_{ij}$-twisted 2-cocycle
equation
\begin{equation}
  \label{eq:coc2}
  \lambda_{ij}(g_{jkl})\, g_{ijl} = g_{ijk} g_{ikl}
\end{equation}
in $G_{i\, \mid  U_{ijkl}}$. 
\end{lemma}
\noindent \noindent {\bf Proof:}   Note that
equation
 \eqref{eq:coc2} is equivalent to the commutativity of the
diagram 
\begin{equation}
\label{def:gijk}
\hspace{1cm} 
  \xymatrix@R=10pt@C=15pt{
   x_i \ar[rr]^{g_{ijl}} \ar[dd]_{g_{ikl}}
    && x_i \ar[dd]^{\lambda_{ij}(g_{jkl})}\\&&\\x_i 
    \ar[rr]_{g_{ijk}}&& x_i
    }\end{equation} 
\noindent above $U_{ijkl}$.
Let us now consider the following 
 cubical diagram: 
\begin{equation}
 \label{verif:2coc}
       \xymatrix@=12pt{
     &&x_l\ar[rrrrr]^{\phi_{jl}} \ar[lldd]_{\phi_{il}}
    \ar '[dd][ddddd]^(.3){\phi_{kl}}
     \ar@{}[rrrrr]_(.2){\,}="1"
     \ar@{}[ddddd]^(.3){\,}="2"
      &&&&& x_j
     \ar[lldd]^{\phi_{ij}}
     \ar@{}[lldd]^(.7){\,}="6"
     \ar[ddddd]^{g_{jkl}}\\
     &&&&&&&\\
     x_i\ar[rrrrr]^(.6){g_{ijl}}
      \ar@{}[rrrrr]^(.7){\,}="5"
     \ar[ddddd]_{g_{ikl}}
     \ar@{}[ddddd]^(.6){\,}="9"
     &&&&&x_i    \ar[ddddd]^(.4){\lambda_{ij}(g_{jkl})}
     \ar@{}[ddddd]_(.8){\,}="17"
     &&\\
     &&&&&&&\\
     &&&&&&&\\
     &&x_k    \ar[lldd]_(.4){\phi_{ik}}
     \ar@{}[lldd]^(.4){\,}="13"
     \ar@{}[lldd]^(.7){\,}="10"
     \ar '[rrr]^(.7){\phi_{jk}}[rrrrr]
     \ar@{}[rrrrr]_(.3){\,}="14"
     &&
     &&&
     x_j    \ar[lldd]^{\phi_{ij}}
     \\
     &&&&&&&
     \\ x_i    \ar[rrrrr]_{g_{ijk}}
      \ar@{}[rrrrr]^(.8){\,}="18"
      &&&&& x_i    &&
     }
\end{equation}
\noindent in which  the left, back, top and bottom
 squares are
of type \eqref{coc0}, and the right-hand one 
of type \eqref{def-lamij}.
Since these
 five        faces are commutative squares, and all the arrows in the 
diagram are
invertible,  the sixth (front) face is also commutative. Since the latter is
simply  the square \eqref{def:gijk}, the lemma is proved.
\begin{flushright}
$\Box$
\end{flushright}
  A pair $(\lambda_{ij}, \, g_{ijk})$ satisfying the equations 
\eqref{eq:coc1}  and \eqref{eq:coc2} :
\begin{equation}
\label{defgerbecoc}
\left \{
 \begin{array}{lc}
  \lambda_{ij} \, \lambda_{jk} &= \ i_{g_{ijk}} \, \lambda_{ik}
   \\
 \lambda_{ij}(g_{jkl})\, g_{ijl} &= \  g_{ijk}\, g_{ikl} \end{array}
\right.
\end{equation}
\noindent is called a $G_i$-valued cocycle pair. 
\noindent It may be viewed as consisting of a
2-cocycle equation for the elements  $g_{ijk}$,
 together with auxiliary data attached to the isomorphisms
$\lambda_{ij}$.
\noindent  However, in contrast  with  the abelian case in which the
inner conjugation term $i_{g_{ijk}}$ is trivial,
 these two equations cannot in general be uncoupled. When such a pair
 is associated to a $G$-gerbe $\pc$ for a fixed bundle of groups $G$,
 the term $\lambda_{ij}$ is a section
 above $U_{ij}$ of the bundle of groups $\mathrm{Aut}_X(G)$, and
 $g_{ijk}$
 is a section of $G$ above $U_{ijk}$. Such pairs
 $(\lambda_{ij},\, g_{ijk})$ will be called $G$-valued cocycle pairs.

\subsection 
\noindent The corresponding coboundary relations
 will now  be worked out by a similar
diagrammatic process. Let us  give ourselves a second
family of local objects $x'_i$  in $\pc_{U_i}$, and of arrows
\begin{equation}
  \label{def:lij1}
\xymatrix@=15pt{
x'_j \ar[rr]^{\phi'_{ij}} && x'_i }
 \end{equation}
above $U_{ij}$. To these correspond by the constructions
\eqref{def-lamij} and \eqref{coc0} a new cocycle pair $(\lambda'_{ij},\,
g'_{ijk})$ satisfying the cocycle relations \eqref{eq:coc1} and \eqref{eq:coc2}.
In order to compare the previous  trivializing data $(x_i, \phi_{ij})$ with
the new one, we also choose (possibly after a harmless  refinement of
the  given
 open cover $\mathcal{U}$ of $X$)  a family of  arrows 
\begin{equation}
\label{def:xi}
\xymatrix@C=30pt{
x_i \ar[r]
^{\chi_i} & x'_i}
\end{equation}
in $\pc_{U_i}$ for all $i$. The lack of compatibility between these
 arrows and the previously chosen arrows \eqref{def:lij} and
 \eqref{def:lij1} is measured by the family of arrows $\te_{ij}:x_i
 \la x_i$ in $\pc_{U_{ij}}$ determined by the commutativity of the
 following diagram:

\begin{equation}
\label{def:deij}
\xymatrix@R=15pt@C=40pt{
x_j \ar[r]^{\phi_{ij}} \ar[dd]_{\chi_j} & x_i \ar[d]^{\chi_i} \\
& x'_i \ar[d]^{\te_{ij}} \\
x'_j \ar[r]_{\phi'_{ij}} & x'_i
}
\end{equation} 
The arrow $\chi_i: x_i \la x'_i$ induces by conjugation an
isomorphism  $r_i:  G_{i}\la G'_{i}$, characterized by the
commutativity of the
 square
\begin{equation}
\label{def-mui}
\xymatrix@R=30pt@C=40pt{ x_i \ar[d]_{\chi_i} \ar[r]^u & x_i 
  \ar[d]^{\chi_i} \\
x'_i  \ar[r]_{r_i(u)} & x'_i
}
\end{equation}
 for all $u \in G_{i}$. The diagram
\eqref{def:deij} therefore conjugates to a diagram
\begin{equation}
\label{def:deij1}
\xymatrix@R=15pt@C=40pt{
G_j \ar[r]^{\lambda_{ij}} \ar[dd]_{r_j} & G_i \ar[d]^{r_i} \\
& G'_i \ar[d]^{i_{\te_{ij}}} \\
G'_j \ar[r]_{\lambda'_{ij}} & G'_i
}
\end{equation} 
above $U_{ij}$ whose commutativity is expressed by the equation
\begin{equation}
\label{eq:cob1}
\lambda'_{ij} = i_{\te_{ij}}\, r_i\, \lambda_{ij} \, r_j^{-1}\,.
\end{equation} 
Consider now the diagram
\begin{equation}
\label{diagcoboun}
\xymatrix@R=15pt@C=50pt{
&x_k \ar[ddl]_{\phi_{jk}} \ar[rr]^{\phi_{ik}}
\ar[ddddd]_(.25){\chi_{k}}   &&x_i \ar[d]^{\chi_i}  \ar[ddl]_{g_{ijk}}\\
&&&  x'_i \ar[ddl]^(.7){r_i(g_{ijk}) } \ar[dddd]^{\te_{ik}} \\
   x_j \ar[rr]_(.65){\phi_{ij}}
\ar[dd]_{\chi_{j}}
&&  x_i
\ar[d]_(.45){\chi_i}
& \\
&& x'_i \ar[d]^{\te_{ij}} &\\
 x'_j  
\ar[rr]^(.65){\phi'_{ij}}
 \ar[ddd]_{\te_{jk}}
 & & x'_i
\ar[ddd]^(.25){ \lambda'_{ij}(\te_{jk})}
 & \\
& x'_k  \ar[ddl]_{\phi'_{jk}} 
\ar[rr]^(.3){\phi'_{ik}} && x'_i \ar[ddl]^{g'_{ijk}} \\
&& &\\
  x'_j \ar[rr]_ {\phi'_{ij}} &&  x'_i & 
}
\end{equation}

\noindent Both the top and the bottom squares commute, since these squares
 are of type  \eqref{coc0}.    So do the back,  the left
 and the top front
vertical squares,
 since all three are of type \eqref{def:deij}. The same is true of
 the lower
  front square, and the upper right vertical square, since these two
 are respectively of the
 form \eqref{def-lamij} and \eqref{def-mui}. It follows that the
 remaining lower right  square in the diagram is also
 commutative, since all the arrows in diagram \eqref{diagcoboun} are 
 invertible.  The commutativity
 of this final square is expressed
 algebraically by the equation
\begin{equation}
\label{eq:cob2}
g'_{ijk} \, \te_{ik}
  = \lambda'_{ij} (\te_{jk})\,\te_{ij} \, r_{i}(g_{ijk}),
\end{equation}
an equation equivalent to  \cite{2-gerbe} (2.4.17).

\bigskip

Let us say that two   cocycle pairs
   $(\lambda_{ij},\, g_{ijk})$ 
and $(\lambda'_{ij}  g'_{ijk})$
 are cohomologous if we are given a pair $(r_i, \te_{ij})$, with
 $r_i \in  \mathrm{Isom}(G_i,\,G'_i)$ and $\te_{ij} \in 
G'_{i\,|U_{ij}}$
satisfying the equations
\begin{equation}
\label{recob}
\left \{
\begin{array}{cl}
\lambda'_{ij}   &=\: i_{\te_{ij}}\,\, r_i\, \lambda_{ij} \, r_j^{-1}\\
g'_{ijk}  &=\:
\lambda'_{ij} (\te_{jk})\,\te_{ij}  \,r_{i}(g_{ijk}) \, \te_{ik}^{-1}
\end{array}
\right.
\end{equation}
\noindent 
Suppose now that $\pc$ is a $G$-gerbe. All the terms in the first
equations in both \eqref{defgerbecoc} and \eqref{recob} are then  elements of 
$\mathrm{Aut} (G)$, while the terms in the corresponding second
equations live in $G$.
The set of equivalence classes of cocycle pairs
\eqref{defgerbecoc}, for the equivalence relation defined by equations
\eqref{recob}, is then denoted $H^1(\mathcal{U},\, G \la
\text{Aut}(G))$,
a notation consistent with that introduced in \S \ref{monoid}
  The limit over the open covers $\mathcal{U}$ is the \v{C}ech
  hypercohomology set $\Check{H}^1(X,\, G \la \text{Aut}(G))$.
 We refer to \cite{2-gerbe} \S 2.6 for the inverse  construction, starting
 from a
 \v{C}ech cocycle
 pair, of the corresponding $G$-gerbe\footnote{In \cite{2-gerbe}
   \S 2.7, we explain how this inverse construction extends to the more
 elaborate context of hypercovers, where  a beautiful interplay
between the \v{C}ech and the descent formalisms arises. This is also
discussed, in more simplicial terms, in \cite{festschrift} \S 6.3-6.6.}.  This  hypercohomology
 set therefore 
 classifies $G$-gerbes on $X$ up to equivalence.

\bigskip

\noindent In geometric terms, this can be understood once we introduce the following
 definition, a categorification of the definition  \eqref{deftorsor} 
 of a $G$-torsor:
\begin{definition}
Let $\gc$ be a monoidal stack on $X$. A left  $\gc$-torsor on $X$  is a
stack $\qc$ on $X$ together with a left action functor
\[ \gc \times \qc \la \qc \]
which is coherently associative and satisfies the unit condition, and for
which the induced functor
\[ \gc \times \qc \la \qc \times \qc \]
defined as in  \eqref{deftorsor}
is an equivalence. In addition,
we require that $\qc$ be locally non-empty.
\end{definition} 
 The 
 following three observations, when put together, 
explain in more global terms why $G$-gerbes 
are classified by the set $H^1(X, G \la \text{Aut}(G))$.
\begin{itemize}
\item To a $G$-gerbe $\pc$ on $X$ is associated its ``bundle of frames''
$\mathcal{E}q(\pc, \text{Tors}(G))$, and the latter is a left torsor
under the monoidal stack  
$\mathcal{E}q(\text{Tors}(G),\, 
\text{Tors}(G))$.
\item By the Morita theorem, this monoidal stack is equivalent to the
monoidal stack $\text{Bitors}(G)$ of $G$-bitorsors on $X$.
\item  The cocycle computations leading up  to \eqref{defbitorcoc} imply 
that the 
  monoidal stack $\text{Bitors}(G)$
 is the stack associated  
to the
monoidal prestack defined by the crossed module $G \stackrel{i}{\la}
\mathrm{Aut}(G)$ \eqref{def:inconj}.
\end{itemize}

\medskip

\begin{remark}
{\rm  For a related discussion of non-abelian
cocycles in a homotopy-theoretic context, see 
J. F. Jardine  \cite{jardine} theorem 13 and
\cite{jardine1} \S 4,
 where a
classification of gerbes equivalent to ours is given, including   the
case  in which 
hypercovers are required.
}\end{remark}

\bigskip

\subsection \noindent {\bf  A topological interpretation of a
  $G$-gerbe (\cite{lbschreier} 4.2)}

\vspace{.3cm}
The context here is that of fibrewise topology, in which all
constructions
are done in the category of spaces above a fixed topological space
$X$.
Let $G$ be   a bundle of  groups above $X$ and $B_XG$ its
classifying space, a space above $X$  whose fiber at a point $x \in X$
is
 the classifying
space $BG_x$ of the group $G_x$.
By construction, $B_XG$ is the geometric realization of the
simplicial space over $X$  whose face and degeneracy operators above
$X$ are defined in the usual fashion (but now in the fibrewise
context)
  starting from  the  multiplication   and  diagonal maps
$G \times_X G \la G$  and $G \la G\times_X G$:

\medskip

\[ \xymatrix@C=40pt{
\ldots  \quad
\ar@<-2ex>[r]   \ar@<2ex>[r]  
 &  G \times_X G \ar@<1ex>[r] \ar[r] \ar@<-1ex>[r]
\ar@/_/[drr] 
 &G \ar[dr]
  \ar@<.5ex>[r]
\ar@<-.5ex>[r] \ar@<-2ex>[l]  \ar@<-3ex>[l]
 & \ar@<-1.5ex>[l]
X  \ar@{=}[d] \\
 & &&X
}\]
We   attach to $G$ 
 the bundle  $\text{Eq}_X (BG)$  of  group-like  topological  monoids
   of self-fiber-homotopy equivalences
 of $B_XG$  over $X$. The fibrewise homotopy fiber of the  evaluation map
\[\text{ev}_{X,\,\ast}: \text{Eq}_X (BG) \la BG\,, \]
which associates to such  an equivalence its value 
 at the distinguished section $ \ast$ of $B_XG$ above $X$, is the 
submonoid $\text{Eq}_{X,\, \ast} (B_XG)$ of pointed fibrewise homotopy 
self-equivalences   of $B_XG$ . The 
latter is fiber homotopy equivalent, by the fibrewise 
functor $\pi_1(-,\, \ast)$, to 
the  bundle of groups $\text{Aut}_X(G)$, whose fiber at a point $x \in
X$ is the group $\text{Aut}(G_x)$.
\noindent This   fibration  sequence   of spaces over $X$
\[ \text{Aut}_X(G)  \la  \text{Eq}_X (B_XG) \la B_XG \]
\noindent  is therefore equivalent to 
 a fibration sequence of topological monoids over $X$, 
 the first two
of which are simply  bundles of   groups on $X$
\begin{equation}
\label{borelc}
 G \stackrel{i}{\la}  \text{Aut}(G)  \la  \text{Eq}_X (BG)  \,.
\end{equation}
 This yields an identification of $\text{Eq}_X (BG)$
with the  fibrewise Borel construction\footnote{Our use here of the
  notation $\times^G$ is meant to be close to the topologists'
  $\times_G$. Algebraic geometers often  denote such a
  $G$-equivariant product  by
  $\wedge^G$, as we did in \eqref{def:P-twist}.}
 $E_XG \times^{G}_X \mathrm{Aut}_X(G)$.
Our discussion in \S \ref{monoid}
asserts that this identification  preserves the
multiplications, so long as the multiplication on the Borel
construction
 is given by
an appropriate iterated  semi-direct product construction, whose
 first non-trivial
stage is defined as in 
 \eqref{semi1}. We refer to \cite{lbschreier} for
a somewhat more detailed discussion of this assertion, and to \cite{d-z} \S 4
for a  related  discussion, in the absolute rather than in the fibrewise
context,  of the corresponding fibration sequence 
\[BG \la B \,\text{Aut}(G) \la  B\text{Eq} (BG) \]
 (or rather to  its generalization in which the classifying space
$BG$ replaced by an arbitrary topological space $Y$).
\noindent This proves:
\begin{proposition}
The simplicial group over $X$ associated to the crossed module
$G \la \text{Aut}(G)$  over $X$ associated to a bundle of groups $G$ 
is a model for the group-like topological  monoid  $\text{Eq}_X (BG)$.
\end{proposition}
For any group $G$,  the set  $H^1(X, G \la \text{Aut}(G))$  of  
 1-cocycle classes describes
 the classes of  fibrations  over $X$ which are locally homotopy equivalent to
 the space  $BG$, and the corresponding assertion when  $G$
is  a bundle of groups on $X$  is also true. We refer to the recent
preprint of J. Wirth and J. Stasheff \cite{wi-sta} for a related discussion of
fiber homotopy equivalence classes of locally homotopy trivial
fibrations, also from a cocyclic point of view. 
\begin{example}{\rm 
Let us sketch here a  modernized  proof of O. Schreier's
cocyclic classification \linebreak (in 1926 !) of
 (non-abelian, non-central)
group extensions \cite{schreier}, which is much less well-known than the
special
case in which the extensions are central.

\medskip

\noindent  Consider a short exact sequence of groups
\eqref{sexs}. Applying the classifying space functor $B$,
 this induces a fibration
\[ BG \la BH \stackrel{\pi}{\la} BK \]
of pointed spaces above $BK$, 
and all the fibers of $\pi$ are homotopically equivalent to $BG$. It follows
that this fibration determines an element in the
 pointed set  $H^1(BK,\, G \la \text{Aut}(G))$. Conversely, such a
 cohomology class determines a fibration $E \stackrel{\phi}{\la} BK$
 above $BK$,
  whose fibers
are  homotopy equivalent to the space  $BG$.  Since both $BG$ and $BK$
have distinguished points, so does $E$.
 Applying the fundamental group functor to
this fibration of pointed spaces determines 
 a short  sequence  of  groups
\[1 \la G \la H \la K \la 1 \,. \]
\begin{flushright}
$\Box$
\end{flushright}
}
\end{example}

{\Large
\section{\bf 2-stacks and 2-gerbes}
}
\vspace{.3cm}

 \subsection 
\noindent  We will  now extend the discussion of section 5
from 1- to 2-categories. A 2-groupoid is defined here as a 2-category
whose 1-arrows are invertible up to a 2-arrow, and whose 2-arrows are
strictly invertible. 

 \begin{definition}
  A fibered 2-category in 2-groupoids above a space $X$ 
consists in a family of 2-groupoids $\cc_U$, for each open set $U$ in $X$,
 together with an inverse image 2-functor
\[ f^{\ast}: \cc_U \la \cc_{U_1} \]
\noindent associated to every inclusion of open sets  $f: U_1 \subset U$
 (which is the identity
 whenever $f = 1_U$), and a natural transfomation 
\[ \phi_{f,g}: (fg)^{\ast} \Longrightarrow g^{\ast}\, f^{\ast} \]
for every pair of composable inclusions
\[ U_2 \stackrel{g}{\hookrightarrow} U_1  \stackrel{f}{\hookrightarrow} U \:.\]
\noindent For each triple of composable inclusions 
\[ U_3 \stackrel{h}{\hookrightarrow}   U_2 \stackrel{g}{ \hookrightarrow}
 U_1  \stackrel{f}{ \hookrightarrow} U, \]
\noindent  we require  a modification
\bigskip
\[
\label{defthi}
 \xymatrix{
                 (fgh)^{\ast}           
\ar@2@/^1.5pc/[rrr]^{\psi_{f,g,h}}_{\,}="1" 
\ar@2@/_1.5pc/[rrr]_{\chi_{f,g,h}}_{\,}="2" &
& &  h^{\ast}    g^{\ast}    f^{\ast}
\ar@{}"1";"2"^(.2){\,}="3" 
\ar@{}"1";"2"^(.8){\,}="4"
\ar@3"3";"4"^{\al_{f,g,h}}
}
\]
 \noindent betweeen the composite natural transfomations
\[ \psi_{f,g,h}: (fgh)^{\ast} \Longrightarrow h^{\ast}\,(fg)^{\ast}
 \Longrightarrow h^{\ast} \, (g^{\ast} f^{\ast}) \] 
\noindent and
\[\chi_{f,g,h}:  (fgh)^{\ast} \Longrightarrow (gh)^{\ast}\, f^{\ast}
 \Longrightarrow (h^{\ast} g^{\ast}) \, f^{\ast}. \]
\noindent Finally, for any  $U_4 \stackrel{k}{\hookrightarrow} U_3$, 
 the two methods  by which the induced modifications $\alpha$ 
compare  the composite
 2-arrows  
\[ (fghk)^{\ast} \Longrightarrow (ghk)^{\ast} f^{\ast} \Longrightarrow 
((hk)^{\ast}g^{\ast}f^{\ast} \Longrightarrow
 k^{\ast}  h^{\ast}  g^{\ast}  f^{\ast} \]
\noindent and 
\[  (fghk)^{\ast}  \Longrightarrow k^{\ast} (fgh)^{\ast} \Longrightarrow 
k^{\ast} (h^{\ast}(fg)^{\ast}) \Longrightarrow 
k^{\ast}  h^{\ast}  g^{\ast}  f^{\ast} \]
\noindent   must  coincide.
\end{definition}

\begin{definition}
  A 2-stack in 2-groupoids  above a space $X$ is a fibered 2-category
 in 2-groupoids above $X$ such that
\begin{itemize}
 \item For every pair of objects $X, Y 
\in \cc_U$, the fibered category  $ \mathrm{Ar}_{\cc_U}(X,\,Y)$ 
is a stack on $U$.

\item 2-descent is effective  for objects in  $\cc$. 
 \end{itemize}
\end{definition}
\noindent The 2-descent condition asserts that  we are given, for  an open
 covering  $(U_{\alpha})_{\alpha \in J}$
 of an open set $U \subset X$,  a 
family of objects $x_{\alpha} \in \cc_{U_{\alpha}}$,
 of 1-arrows $ \phi_{\alpha \beta}:x_\alpha \la x_{\beta}$
 between the restrictions
 to $\cc_{U_{\alpha \beta}}$ of the objects $x_{\alpha}$ and
  $x_{\beta }$, and a family of 2-arrow 
\begin{equation}
\label{tetra}
\xymatrix@C=25pt@R=20pt{
&x_{\be} \ar[dl]_{\phi_{\be \ga}}
\ar[dr]^{\phi_{\al \be}}
 \ar@{}[d]^(.5){\,}="1" \ar@{}[d]^(.9){\,}="3"
 &\\
x_{\ga} \ar[rr]_{\phi_
{\al \ga}}
  \ar@{}[r]^(.99){\,}="2"   
 &&x_{\alpha}
\ar@{=>}"1";"3"_{ \psi_{\al \be \ga}}
}\end{equation}

\vspace{.25cm}

\noindent for which the 
   tetrahedral diagram of 2-arrows  whose four faces are the
   restrictions of the requisite
   2-arrows  $\psi$ 
\eqref{tetra}  to  $\cc_{U_{\al \be \ga \de}}$ commutes:

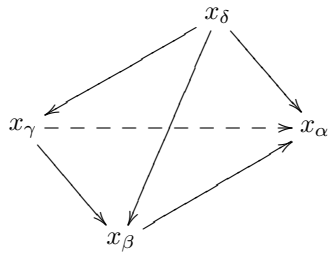
\vspace{.25cm}

\[ \xymatrix@R=30pt@C=20pt{&&x_{\delta} \ar[dll]
 \ar[ddl] \ar[dr] &\\
x_{\gamma} \ar@{-->}[rrr] \ar[rd] &&&x_{\al} \\
& x_{\be} \ar[rru] &&
}\]

\noindent The 2-descent condition
  $(x_\alpha,\, \phi_{\al \be},\,\psi_{\al \be \ga})$ is effective
 if there exists an object $x \in \cc_U$, together with  1-arrows
 $x_{|U_{\al}} \simeq x_{\al}$ in $\cc_{U_\alpha}$ which are
  compatible with the given 1- and
 2-arrows  $\phi_{\al,\be}$ and
  $\psi_{\al \be \ga}$.

\begin{definition}
A 2-gerbe $\mathcal{P}$  is a 2-stack  in 2-groupoids on $X$ 
which is locally non-empty and locally connected. 

\end{definition}

\vspace{.3cm}

\noindent To each object $x$ in $\pc_{U} $ is associated a  group like 
 monoidal stack (or $gr$-stack)  $\gc_x:= \mathcal{A}r_{U}(x,\,x)$
 above $U$.

\vspace{.5cm}

\begin{definition}
Let $\gc$ be a group-like monoidal stack on $X$. 
We say that  a 2-gerbe $\pc$ is a $\gc$-2-gerbe if there exists an open covering
  $\mathcal{U}:=(U_i)_{i \in I}$ of $X$, a family of objects $x_i \in
 \pc_{U_i} $, and $U_i$-equivalences $\gc_{U_i} \simeq \gc_{x_i}$.
\end{definition}

\bigskip

\subsection
\noindent {\bf Cocycles for 2-gerbes
:}

\bigskip

In order to obtain a cocyclic description of a $\gc$-2-gerbe $\pc$, we will
now categorify the constructions in \S 5.
 We choose  paths 
\begin{equation}
\label{2phij}
\phi_{ij}:x_j \la x_i
\end{equation}
  in the 2-groupoid 
 $\pc_{U_{ij}}$, together with  quasi-inverses $x_i \la x_j$ and
 pairs of 2-arrows
\begin{equation}
\label{def:rs}
\begin{array}{ccc}
\xymatrix{\phi_{ij}\, \phi_{ij}^{-1} \ar@{=>}[r]^(.6){r_{ij}} & 1_{x_i}} 
&\hspace{2cm} & 
\xymatrix{ \phi_{ij}^{-1}\, \phi_{ij} \ar@{=>}[r]^(.6){s_{ij}} & 1_{x_j}} \,.
\end{array}
\end{equation}
 These 
 determine  a monoidal equivalence
\begin{equation}
\label{def-lamij0} 
\lambda_{ij}: 
\gc_{|U_{ij}} \la \gc_{|U_{ij}}
\end{equation}
 as well as, functorially  each object 
$\ga \in \gc_{|U_{ij}}$, a 2-arrow $M_{ij}(\ga)$
\begin{equation}
        \label{2def-lamij}
        \xymatrix@R10pt@C=15pt{
   x_j \ar[rr]^{\ga} \ar[dd]_{\phi_{ij}}
    && x_j \ar[dd]^{\phi_{ij}}_{}="2"\\&&\\x_i 
    \ar[rr]_{\lambda_{ij}(\ga)}^{}="1"&& x_i
\ar@{}"1";"2"^(.2){}="4"
\ar@{}"1";"2"^(.85){}="5"
\ar@{=>}"4";"5"
    }\end{equation}
which categorifies diagram \eqref{def-lamij}.
 In fact, the 2-arrows $r$ and $s$ \eqref{def:rs}
can be chosen coherently, and the induced 2-arrow
\eqref{2def-lamij} therefore  does not carry any additional  cohomological
information. For this reason,we will not label such a  2-arrow
$M_{ij}(\ga)$
explicitly when it  occurs
in one of our diagrams. For the same reason, 
 we  will  
 treat diagrams such as  \eqref{2def-lamij}  as
commutative squares.

\bigskip

\noindent The paths $\phi_{ij}$  and their inverses also give us
 objects $g_{ijk} \in \gc_{U_{ijk}}$ and   2-arrows $m_{ijk}$:
 \begin{equation}
         \label{def:mijk}
 \hspace{-.9cm}        \xymatrix@R=12pt@C=20pt{
    x_k  \ar[rr]^{\phi_{jk}} \ar[dd]_{\phi_{ik}}
     &&x_j \ar[dd]^{\phi_{ij}}
     \ar@{}[dd]_(.55){\, }="1"\\&&\\ x_i \ar[rr]_{g_{ijk}}
     \ar@{}[rr]^(.55){\,}="2"&& x_i 
      \ar@{}"1";"2"^(.2){\,}="3"
      \ar@{}"1";"2"^(.8){\,}="4"
      \ar@{=>}"4";"3"^{m_{ijk}}
     }\end{equation}
\noindent These in turn determine a 2-arrow $\nu_{ijkl}$ above $U_{ijkl}$ 
 \begin{equation}
         \label{def:nuijkl}
         \xymatrix@R=12pt@C=20pt{
    x_i \ar[rr]^{g_{ijl}} \ar[dd]_{g_{ikl}}
     && x_i \ar[dd]^{\lambda_{ij}(g_{jkl})}
     \ar@{}[dd]_(.55){\, }="1"\\&&\\ x_i \ar[rr]_{g_{ijk}}
     \ar@{}[rr]^(.55){\,}="2"&& x_i 
      \ar@{}"1";"2"^(.2){\,}="3"
      \ar@{}"1";"2"^(.8){\,}="4"
      \ar@{=>}"3";"4"_{\nu_{ijkl}}
     }
  \end{equation}
\noindent as the unique
2-arrow such that the following diagram of 2-arrows with right-hand
face \eqref{2def-lamij} and front face
 $\nu_{ijkl}$ commutes:
\begin{equation}
   \label{imacube}
         \xymatrix@=14pt{
    && x_l  \ar[rrrrr]^{\phi_{jl}} \ar[lldd]_{\phi_{il}}
   \ar '[dd][ddddd]^(.3){\phi_{kl}}
    \ar@{}[rrrrr]_(.2){\,}="1"
    \ar@{}[ddddd]^(.3){\,}="2"
     &&&&&x_j
    \ar[lldd]_(.6){\phi_{ij}}
    \ar@{}[lldd]^(.7){\,}="6"
 \ar@{}[lldd]^(.6){\,}="21"
 \ar[ddddd]^{g_{jkl}}
 \ar@{}[ddddd]^(.6){\,}="22"\\
    &&&&&&&\\
   x_i  \ar[rrrrr]^(.6){g_{ijl}}
     \ar@{}[rrrrr]^(.7){\,}="5"
    \ar[ddddd]_{g_{ikl}}
    \ar@{}[ddddd]^(.6){\,}="9"
    &&&&&x_i \ar[ddddd]^(.4){\lambda_{ij}(g_{jkl})} 
    \ar@{}[ddddd]_(.8){\,}="17"
    &&\\
    &&&&&&&\\
    &&&&&&&\\
    && x_k \ar[lldd]^(.4){\phi_{ik}} 
    \ar@{}[lldd]^(.4){\,}="13"
    \ar@{}[lldd]^(.7){\,}="10"
    \ar '[rrr]^(.7){\phi_{jk}} [rrrrr]
    \ar@{}[rrrrr]_(.3){\,}="14"
    &&
    &&&
    x_j \ar[lldd]^{\phi_{ij}}
    \\
    &&&&&&&
    \\ x_i  \ar[rrrrr]_{g_{ijk}}
     \ar@{}[rrrrr]^(.8){\,}="18"
     &&&&& x_i && 
    \ar@{}"1";"2"^(.3){\,}="3"
     \ar@{}"1";"2"^(.8){\,}="4"
   \ar@{=>}"3";"4"^{m_{jkl}}
   \ar@{}"5";"6"^(.4){\,}="7"
   \ar@{}"5";"6"^(.75){\,}="8"
   \ar@{}"9";"10"^(.1){\,}="11"
   \ar@{}"9";"10"^(.5){\,}="12"
   \ar@{=>}"7";"8"^{m_{ijl}}
   \ar@{=>}"11";"12"^{m_{ikl}}
    \ar@{}"13";"14"^(.3){\,}="15"
     \ar@{}"13";"14"^(.6){\,}="16"
    \ar@{=>}"15";"16"_{m_{ijk}}
     \ar@{}"17";"18"^(.2){\,}="19"
       \ar@{}"17";"18"^(.8){\,}="20"
        \ar@{=>}"19";"20"_{\nu_{ijkl}}
        \ar@{}"21";"22"^(.35){\,}="23"
\ar@{}"21";"22"^(.25){\,}="25"
        \ar@{}"21";"22"^(.8){\,}="24"
 \ar@{}"21";"22"^(.7){\,}="26"
  \ar@{}"23";"24"_{}
\ar@{=>}"25";"26"
}
\end{equation}

\noindent This cube  in $\pc_{U_{ijkl}}$
will be denoted  $C_{ijkl}$.
Consider now the following diagram:

   \begin{equation}
 \label{imabianchicube}
         \xymatrix@=15pt{
    &&x_i \ar[rrrrr]^{g_{ijl}} 
   \ar '[dd][ddddd]^(.3){g_{ikl}}
    \ar@{}[ddddd]_(.8){\,}="30"
    \ar@{}[rrrrr]_(.2){\,}="1"
    \ar@{}[ddddd]^(.3){\,}="2"
     &&&&& x_i
    \ar[ddddd]^{\lambda_{ij}(g_{jkl})}\\
    &&&&&&&\\
   x_i \ar[rrrrr]^(.65){g_{ijm}}
    \ar@{}[rrrrr]^(.3){\,}="25"
     \ar[uurr]^{g_{ilm}}
     \ar@{}[uurr]^(.3){\,}="26"
     \ar@{}[rrrrr]^(.7){\,}="5"
    \ar[ddddd]_{g_{ikm}}
    \ar@{}[ddddd]^(.6){\,}="9"
    &&&&&x_i
    \ar[uurr]^{\lambda_{ij}(g_{jlm})}
     \ar@{}[uurr]_(.3){\,}="22"
      \ar[ddddd]^(.4){\lambda_{ij}(g_{jkm})}
    \ar@{}[ddddd]^(.3){\,}="21"
    \ar@{}[ddddd]_(.8){\,}="17"
    &&\\
    &&&&&&&\\
    &&&&&&&\\
    && x_i
    \ar '[rrr]^(.7){g_{ijk}} [rrrrr]
    \ar@{}[rrrrr]_(.3){\,}="14"
    &&
    &&&
   x_i
    \\
    &&&&&&&
    \\ x_i
    \ar[uurr]_(.8){\lambda_{ik}(g_{klm})}
     \ar@{}[uurr]_(.4){\,}="34"
      \ar@{}[uurr]^(.7){\,}="29"
    \ar[rrrrr]_{g_{ijk}}
     \ar@{}[rrrrr]^(.3){\,}="33"
      \ar@{}[rrrrr]^(.8){\,}="18"
     &&&&&x_i
      \ar[uurr]_(.6){\lambda_{ij}\lambda_{jk}(g_{klm})}
      && 
    \ar@{}"1";"2"^(.3){\,}="3"
     \ar@{}"1";"2"^(.8){\,}="4"
   \ar@{=>}"3";"4"^{\nu_{ijkl}}
   \ar@{}"25";"26"^(.2){\,}="27"
   \ar@{}"25";"26"^(.8){\,}="28"
   \ar@{=>}"27";"28"_{\nu_{ijlm}}
   \ar@{}"29";"30"^(.3){\,}="31"
   \ar@{}"29";"30"^(.7){\,}="32"
    \ar@{=>}"31";"32"^{\nu_{iklm}}
    \ar@{}"21";"22"^(.3){\,}="23"
    \ar@{}"21";"22"^(.7){\,}="24"
    \ar@{=>}"23";"24"_{\lambda_{ij}(\nu_{jklm})}
     \ar@{}"17";"18"^(.2){\,}="19"
       \ar@{}"17";"18"^(.8){\,}="20"
        \ar@{=>}"19";"20"_{\nu_{ijkm}}
         \ar@{}"33";"34"^(.05){\,}="35"
         \ar@{}"33";"34"^(.95){\,}="36"
     \ar@{=>}"35";"36"_{\{\widetilde{m}_{ijk},\,g_{klm}\}^{-1}}
        } 
\end{equation}

\vspace{1.5cm}

\noindent In order to avoid any possible ambiguity, we spell out in
the following table 
the names of the faces of the cube \eqref{imabianchicube}: 

\pagebreak

  \begin{table}[ht]
\label{imanutable}
      \begin{center}
         \renewcommand{\arraystretch}{1.3}
         \begin{tabular}{|c|c|c|c|c|c|}\hline
       left & right &top & bottom & front & back\\
             \hline
             $\nu_{iklm}$ & $\lambda_{ij}(\nu_{jklm})$ & $\nu_{ijlm}$
 &  $\{\widetilde{m}_{ijk},\, g_{klm}\}^{-1}$
             & $\nu_{ijkm}$ & $\nu_{ijkl}$   \\ \hline
       \end{tabular}
\vspace{.2cm}
\caption{The faces of cube   \eqref{imabianchicube}}
\end{center}
\end{table}

\bigskip

\noindent As we see from   this table,   five
 of its faces are defined by arrows $\nu$ \eqref{def:nuijkl}. The
remaining
 bottom 2-arrow $\{\widetilde{m}_{ijk},\, g_{klm}\}^{-1}$ is
 essentially the
inverse of the 
2-arrow $\widetilde{m}_{ijk}(g_{klm})$ 
obtained by evaluating the natural transformation 
\begin{equation}
\label{def:mtilde0}
\widetilde{m}_{ijk}: i_{g_{ijk}}\lambda_{ik} \Rightarrow \lambda_{ij} \, 
\lambda_{jk}
\end{equation}
\noindent induced by conjugation from the 2-arrow $m_{ijk}$  \eqref{def:mijk}
on the object $g_{klm} \in G:=  \text{Aut}_{\pc}(x_k)$. More
precisely, if we compose the latter 2-arrow as follows with the
unlabelled 2-arrow  $M_{g_{ijk}}(\lam_{ik}(g_{klm}))$ associated to $i_{g_{ijk}}$:
\begin{equation}
             \label{def:mtilde}
     \xymatrix@R=50pt@C=65pt{x_i
\ar@{}[rrd]^(.3){\,}="9"
\ar@{}[rrd]^(.45){\,}="10"
 \ar[rr]^{g_{ijk}}&& x_i \\
x_i \ar[u]^{\lambda_{ik}(g_{klm})}
\ar@{}[u]_(.48){\,}="1"
     \ar[rr]_{g_{ijk}}  &&x_i
     \ar@/^2.7pc/[u]^(.4){i_{g_{ijk}}\lambda_{ik}(g_{klm})}
     \ar@{}@/^2.7pc/[u]_{\,}="5"
     \ar@{}@/^2.7pc/[u]^(.48){\,}="2"
     \ar@/_2.7pc/[u]_{\lambda_{ij}\lambda_{jk}(g_{klm})}
     \ar@{}@/_2.7pc/[u]^{\,}="6"
     \ar@{}"1";"2"^(.3){\,}="3"
      \ar@{}"1";"2"^(.4){\,}="4"
      \ar@{}"5";"6"^(.1){\,}="7"
  \ar@{}"5";"6"^(.9){\,}="8"
     \ar@{=>}"7";"8"^{\widetilde{m}_{ijk}(g_{klm})}  
\ar@{=>}"9";"10"^{}
} 
 \end{equation}
we obtain a 2-arrow
\begin{equation}
\label{def:brack}
\xymatrix@=40pt{
x_i \ar@{}[r]_(.1){\,}="1"
 \ar[r]^{g_{ijk}} & x_i \\
x_i\ar[u]^{\lambda_{ik}(g_{klm})} \ar[r]_{g_{ijk}} & x_i\ar[u]_{\lambda_{ij}
\lambda_{ik}(g_{klm})} \ar@{}[u]^(.2){\,}="2"
\ar@{}"1";"2"^(.35){\,}="3"
 \ar@{}"1";"2"^(.75){\,}="4"
 \ar@{=>}"3";"4"
}
\end{equation}
 which we denote by $\{\widetilde{m}_{ijk},\,g_{klm}\}$. 
 It may be  characterized as the unique 2-arrow such that the
 cube
\begin{equation} 
 \label{imabianchicube-bracket}
   \xymatrix@=15pt{
    &&x_k\ar[rrrrr]^{\phi_{jk}} 
   \ar '[dd][ddddd]^(.3){\phi_{ik}}
 \ar@{}[rrrrr]_(.3){\,}="45"
    \ar@{}[ddddd]^(.3){\,}="46"
     &&&&& x_j
\ar[ddddd]^{\phi_{ij}}
\\
    &&&&&&&\\
    x_k\ar[rrrrr]^(.65){\phi_{jk}}
 \ar@{}[rrrrr]_(.3){\,}="50"
    \ar@{}[ddddd]^(.3){\,}="49"
     \ar[uurr]^{g_{klm}}
\ar[ddddd]_{\phi_{ik}}
    &&&&&x_j
 \ar[uurr]^{\lambda_{jk}(g_{klm})}
    \ar[ddddd]^(.4){\phi_{ij}}
    &&\\
    &&&&&&&\\
    &&&&&&&\\
    &&x_i
    \ar '[rrr]^(.7){g_{ijk}} [rrrrr]
    \ar@{}[rrrrr]_(.3){\,}="14"
    &&
    &&&
    x_i
    \\
    &&&&&&&
    \\ x_i 
    \ar[uurr]_(.75){\lambda_{ik}(g_{klm})}
      \ar@{}[uurr]^(.3){\,}="42"
    \ar[rrrrr]_{g_{ijk}}
     \ar@{}[rrrrr]^(.4){\,}="41"
     &&&&& x_i
      \ar[uurr]_(.6){\lambda_{ij}\lambda_{jk}(g_{klm})}
      &&
\ar@{}"41";"42"^(.2){\,}="43"
 \ar@{}"41";"42"^(.7){\,}="44"
 \ar@{=>}"44";"43"^{\;\;\;\;\;\{\widetilde{m}_{ijk},\,g_{jkl}\}}
\ar@{}"45";"46"^(.2){\,}="47"
 \ar@{}"45";"46"^(.7){\,}="48"
 \ar@{=>}"48";"47"_{m_{ijk}}
\ar@{}"49";"50"^(.2){\,}="51"
 \ar@{}"49";"50"^(.7){\,}="52"
 \ar@{=>}"51";"52"_{m_{ijk}}
      }
\end{equation}
 (with three unlabelled faces of type \eqref{2def-lamij}) is commutative.
 For that reason, this  cube
will be denoted
  $\{\,,\,\}$.
 The following proposition provides a geometric
interpretation for the cocycle equation which the 2-arrows
$\nu_{ijkl}$ satisfy.

\begin{proposition}
\label{2com}
The diagram of 2-arrows \eqref{imabianchicube} is commutative.
\end{proposition}

{\large
{\bf Proof:}
}
Consider  the following hypercubic
diagram, from which the 2-arrows have all been omitted for greater
legibility.

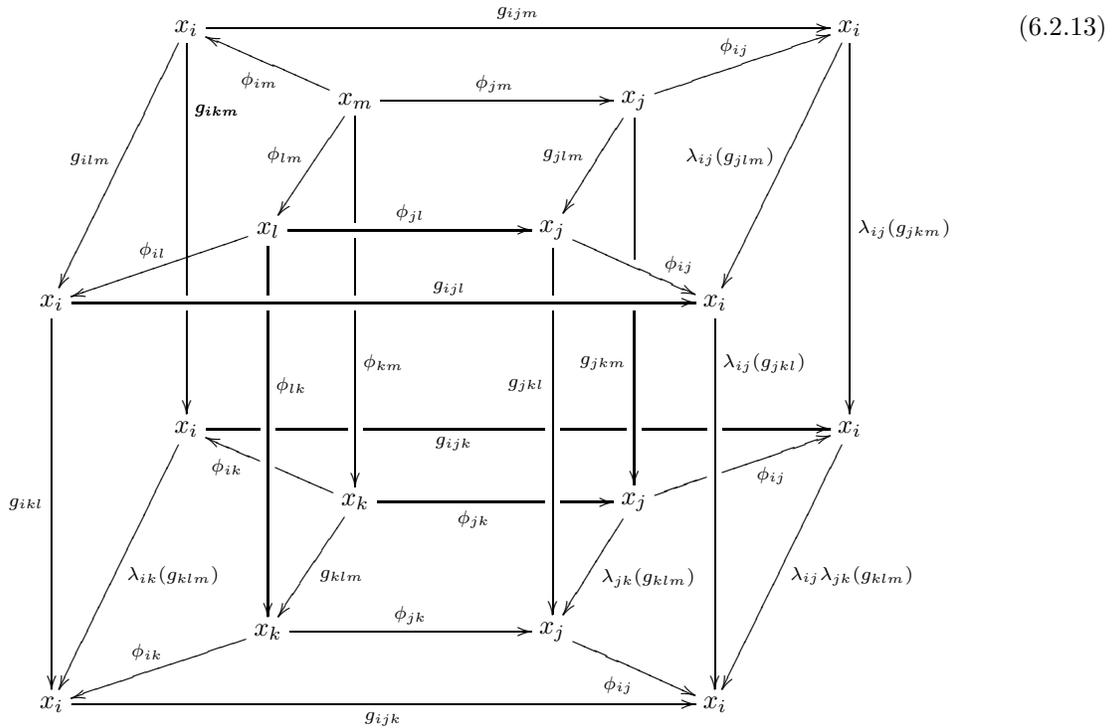
\vspace{1.5cm}

    \begin{equation}  \label{imahypercube}
\xymatrix@=15pt{
   &&x_i  \ar[lldddd]_{g_{ilm}}
   \ar[rrrrrrrrr]^{g_{ijm}}
 \ar@{-}[dddd]^(.3){g_{ikm}} |!{[ddddll];[dddr]}\hole
 \ar@{-}[dddd]^(.3){g_{ikm}}
   \ar@{<-}[drr]_{\phi_{im}}
   &&&&&&&&&x_i \ar[lldddd]_{\lambda_{ij}(g_{jlm})}
   \ar[dddddd]^{\lambda_{ij}(g_{jkm})}
   \ar@{<-}[dlll]_{\phi_{ij}}
   \\&&&& x_m \ar[ldd]_{\phi_{lm}} 
   \ar'[dd]'[ddd][dddddd]^(.3){\phi_{km}}
   \ar[rrrr]^{\phi_{jm}} &&&& x_j 
   \ar[ddl]_{g_{jlm}}
   \ar@{-}[ddd] |!{[ddl];[dddr]}\hole  
   &&&
   \\ &&&  &&& &&& &&
   \\&&&x_l 
   \ar[rrrr]^{\phi_{jl}}\ar@{-}[d]
   &&&& x_j
   \ar'[d][dddddd]_(.25){g_{jkl}} &&&&
   \\
    x_i\ar[rrrrrrrrr]^(.6){g_{ijl}}
   \ar[dddddd]_{g_{ikl}} \ar@{<-}[urrr]^{\phi_{il}}
   && \ar[dd] \ar@{-}[u] |!{[ll];[ur]}\hole
   &\ar[ddddd]^(.25){\phi_{lk}} &&& &&
   \ar[ddd]_(.3){g_{jkm}}
   &  x_i \ar[dddddd]^(.15){\lambda_{ij}(g_{jkl})}
   \ar@{<-}[llu]_(.3){\,\phi_{ij}}&&
   \\
   &&&&&&&&&&&
   \\&&x_i \ar@{<-}[drr]_(.3){\phi_{ik}} |!{[uuur];[dddr]}\hole
   \ar[ddddll]^{\lambda_{ik}(g_{klm})}
   \ar'[r]'[rr]'[rrrrr]_{g_{ijk}}'[rrrrrr]'[rrrrrrr][rrrrrrrrr]
   &&& &&& &&&x_i
   \ar[ddddll]^{\lambda_{ij}\lambda_{jk}(g_{klm})} 
   \ar@{<-}[dlll]^(.4){\phi_{ij}} |!{[lldd];[lluuuu]}\hole
   \\&&&&x_k \ar[ddl]^{g_{klm}}
   \ar'[rrr]_(.6){\phi_{jk}}[rrrr]
   &&&& x_j \ar[ddl]^{\lambda_{jk}(g_{klm})} &&& 
   \\&&&&&&&&&&&
   \\&&&x_k  \ar[rrrr]^{\phi_{jk}} &&& &x_j&&&
   \\
  x_i  \ar@{<-}[rrru]^{\phi_{ik}\,\,}
   \ar[rrrrrrrrr]_{g_{ijk}} &&&&&&&&& x_i 
   \ar@{<-}[llu]^{\phi_{ij}\,\,}&& 
   }
\end{equation}

\vspace{1.5cm}

\noindent The following table is provided as a help  in understanding 
diagram \eqref{imahypercube}.
 The first line describes the position in the hypercube of
  each of the eight cubes from which
it has been constructed, and the middle line gives  each of these a name.
Finally, the last line describes the face by which it is attached to
the inner cube $C_{jklm}$.

   \begin{table}[ht]
\label{imahypertable}
      \begin{center}
         \renewcommand{\arraystretch}{1.3}
         \begin{tabular}{|c|c|c|c|c|c|c|c|}\hline
           inner  & left & right &top & bottom & front & back& outer\\
             \hline
             $C_{jklm}$ & $C_{iklm}$ &$\text{Conj}(\phi_{ij})$ & $C_{ijlm}$
             & $\{\:,\:\}$ & $C_{ijkl}$ & $C_{ijkm}$ &
             \eqref{imabianchicube} 
 \\ \hline
              &$m_{klm}$ &$\nu_{jklm}$
               &$m_{jlm}$&$M_{jk}(m_{klm})$& $m_{jkl}$&$m_{jkm}$&
             \\ \hline
       \end{tabular}
\vspace{.2cm}
\caption{The constituent cubes of diagram \eqref{imahypercube}}
\end{center}
\end{table}

\bigskip

\noindent  Only one cube in this table has not yet been
described. It  is the cube $\mathrm{Conj}(\phi_{ij})$
 which appears on the right in
diagram \eqref{imahypercube}. It  describes the construction of 
 the 2-arrow
 $\lambda_{ij}(\nu_{jklm})$  starting from
  $\nu_{jklm}$, by
conjugation  of its source and target arrows by  the  1-arrows  $\phi_{ij}$.

\bigskip

 Now that diagram \eqref{imahypercube}  has been properly described, the
proof of  proposition \ref{2com}  is immediate, and goes along the same lines
as the proof of lemma \ref{lem:coc2}. One simply observes that each of the
 first 
seven cubes in table 2 is a commutative diagram of 2-arrows. Since all
their
  constituent  2-arrows are invertible, the remaining outer cube is also a
commutative diagram of 2-arrows.  The latter cube is simply
 \eqref{imabianchicube}, though with a different orientation,   
 so the proof of the proposition is now complete.
\begin{flushright}
$\Box$
\end{flushright}
\begin{remark} {\rm
When $i=j$, it is natural to choose as arrow $\phi_{ij}$
\eqref{2phij}
 the identity arrow $ 1_{x_i}$.  When $i=j$ or $j=k$,
 it is then possible to set
 $g_{ijk} = 1_{x_i}$
 and to choose the identity 2-arrow for
$m_{ijk}$.
These choices  yield  the following normalization conditions:
\[
\begin{cases}
\lambda_{ij} = 1\  \text{whenever} \ i=j \\
g_{ijk}= 1 \ \text{and} \ \widetilde{m}_{ijk} = 1 \ \text{whenever} \
i=j \ \text{or} \ j=k \\
\nu_{ijkl} = 1 \ \text{whenever} \ i=j,\ j= k,\, \text{or} \ k=l
\end{cases}
\]
}
\end{remark}
\subsection
\noindent {\bf Algebraic description of the cocycle conditions:}

\bigskip
In order to obtain a genuinely cocyclic description
of a $\gc$-2-gerbe, it is necessary to translate  
 proposition  \ref{2com} into an algebraic statement. 
As a preliminary step, we implement such a translation  for the cubical  diagram
$C_{ijkl}$ \eqref{imacube} by which we defined the 2-arrow
$\nu_{ijkl}$. We reproduce this cube as 
\begin{equation}
   \label{imacube-dot}
         \xymatrix@=14pt{
    && *+[F]{x_{l}} \ar[rrrrr]^{\phi_{jl}} \ar@{-->}[lldd]_{\phi_{il}}
   \ar@2{.>} '[dd][ddddd]^(.3){\phi_{kl}}
    \ar@{}[rrrrr]_(.2){\,}="1"
    \ar@{}[ddddd]^(.3){\,}="2"
     &&&&&x_j
    \ar[lldd]_(.6){\phi_{ij}}
    \ar@{}[lldd]^(.7){\,}="6"
 \ar@{}[lldd]^(.6){\,}="21"
 \ar[ddddd]^{g_{jkl}}
 \ar@{}[ddddd]^(.6){\,}="22"\\
    &&&&&&&\\
   x_i  \ar@{-->}[rrrrr]^(.6){g_{ijl}}
     \ar@{}[rrrrr]^(.7){\,}="5"
    \ar[ddddd]_{g_{ikl}}
    \ar@{}[ddddd]^(.6){\,}="9"
    &&&&&x_i \ar@{-->}[ddddd]^(.4){\lambda_{ij}(g_{jkl})} 
    \ar@{}[ddddd]_(.8){\,}="17"
    &&\\
    &&&&&&&\\
    &&&&&&&\\
    && x_k \ar[lldd]^(.4){\phi_{ik}} 
    \ar@{}[lldd]^(.4){\,}="13"
    \ar@{}[lldd]^(.7){\,}="10"
    \ar@2{.>} '[rrr]^(.7){\phi_{jk}} [rrrrr]
    \ar@{}[rrrrr]_(.3){\,}="14"
    &&
    &&&
    x_j \ar@2{.>}[lldd]^{\phi_{ij}}
    \\
    &&&&&&&
    \\x_i  \ar[rrrrr]_{g_{ijk}}
     \ar@{}[rrrrr]^(.8){\,}="18"
     &&&&&  *+[F]{x_i} && 
    \ar@{}"1";"2"^(.3){\,}="3"
     \ar@{}"1";"2"^(.8){\,}="4"
   \ar@{=>}"3";"4"^{m_{jkl}}
   \ar@{}"5";"6"^(.4){\,}="7"
   \ar@{}"5";"6"^(.75){\,}="8"
   \ar@{}"9";"10"^(.1){\,}="11"
   \ar@{}"9";"10"^(.5){\,}="12"
   \ar@{=>}"7";"8"^{m_{ijl}}
   \ar@{=>}"11";"12"^{m_{ikl}}
    \ar@{}"13";"14"^(.3){\,}="15"
     \ar@{}"13";"14"^(.6){\,}="16"
    \ar@{=>}"15";"16"_{m_{ijk}}
     \ar@{}"17";"18"^(.2){\,}="19"
       \ar@{}"17";"18"^(.8){\,}="20"
        \ar@{=>}"19";"20"_{\nu_{ijkl}}
        \ar@{}"21";"22"^(.35){\,}="23"
\ar@{}"21";"22"^(.25){\,}="25"
        \ar@{}"21";"22"^(.8){\,}="24"
 \ar@{}"21";"22"^(.7){\,}="26"
  \ar@{}"23";"24"_{}
\ar@{=>}"25";"26"
}
\end{equation}
and  consider the two  
 composite paths of 1-arrows from the framed vertex $x_l$ to the 
framed vertex $x_i$ respectively displayed by arrows of type
$\xymatrix{\ar@{-->}[r]&}$ and  $\xymatrix{\ar@2{.>}[r]&}$.  The
commutativity of our cube is equivalent to  the assertion that the two
possible composite 2-arrows  from the path $\xymatrix{\ar@{-->}[r]&}$  to the
path $\xymatrix{\ar@2{.>}[r]&}$
coincide. This assertion translates, when taking into
 account the whiskerings which arise
whenever one  considers a face of the cube which does not
 contain the framed vertex $x_i$,
to the  equation
\begin{equation}
\label{defnualg}
 m_{ijk} \, (g_{ijk} \ast m_{ikl})\, \nu_{ijkl} =
(\phi_{ij} \ast m_{jkl})\, (\lam_{ij}(g_{jkl}) \ast m_{ijl})
\end{equation}
which  algebraically  defines  the 2-arrow $\nu_{ijkl}$ in terms of the
2-arrows of type $m_{ijk}$ \eqref{def:mijk}.
For reasons which will appear later on, we have neglected here the
whiskerings by 1-arrows on the right, for faces of the cube which do
not contain the framed vertex $x_l$ from which all paths considered originate. 
With the left-hand side of this equality labelled ``1'' and the
right-hand side ``2'', the two sides are compared according
 to the following scheme in the 2-category $\pc_{U_{ijkl}}$ : 
\[\xymatrix@=70pt{
*+[F]{x_l} \ar@/^2pc/[rr]^{\lambda_{ij}(g_{jkl})\,g_{ijl}\,\phi_{il}} 
 \ar@/_2pc/[rr]_{\phi_{ij}\, \phi_{jk}\, \phi_{jl}} 
\ar@/^2pc/[rr]^(.3){}="1" \ar@/^2pc/[rr]^(.7){}="2"
  \ar@/_2pc/[rr]^(.3){}="3"\ar@/_2pc/[rr]^(.7){}="4"
 && *+[F]{x_i}
\ar@{}"1";"3"^(.2){\,}="5"
\ar@{}"1";"3"^(.8){\,}="6"
\ar@{=>}"5";"6"_{1}
\ar@{}"2";"4"^(.2){\,}="7"
\ar@{}"2";"4"^(.8){\,}="8"
\ar@{=>}"7";"8"^{2}
\ar@{}"5";"6"^(.5){\,}="9"
\ar@{}"7";"8"_(.5){\,}="10"
\ar@{}"9";"10"^(.3){\,}="11"
\ar@{}"9";"10"^(.7){\,}="12"
\ar@{=}"11";"12"
}\]

\noindent Consider now  a 2-arrow
\begin{equation}
\label{def:m}
\xymatrix{
y \ar@/^1pc/[rr]^{\al} \ar@/_1pc/[rr]_{\be}
\ar@/^1pc/[rr]_(.45){\,}="1" 
\ar@/_1pc/[rr]_(.45){\,}="2"
\ar@{}"1";"2"^(.1){\,}="3"
\ar@{}"1";"2"^(.8){\,}="4"
\ar@{=>}"3";"4"_{^m}
&& x
}
\end{equation}
in $\pc_U$, and denote by $\alpha_\ast$ and $\beta_\ast$ 
the functors $\gc_U \la \gc_U$ which conjugation by  $\al$
 and $\be$ respectively
define. The conjugate of  any  1-arrow  $u \in
\mathrm{ob} \,\gc_U = \mathrm{Ar}_{\pc_U}(y,\,y)$ by
 the 2-arrow $m$ is the composite  2-arrow 
\begin{equation}
\label{def:m1}
\xymatrix@=40pt{
x \ar@/^1pc/[rr]^{\al^{-1}} \ar@/_1pc/[rr]_{\be^{-1}}
\ar@/^1pc/[rr]_(.45){\,}="5" 
\ar@/_1pc/[rr]_(.45){\,}="6"
\ar@{}"5";"6"^(.1){\,}="7"
\ar@{}"5";"6"^(.8){\,}="8"
\ar@{=>}"7";"8"_{m^{-1}}
&& y \ar[r]^u & y 
 \ar@/^1pc/[rr]^{\al} \ar@/_1pc/[rr]_{\be}
\ar@/^1pc/[rr]_(.45){\,}="1" 
\ar@/_1pc/[rr]_(.45){\,}="2"
\ar@{}"1";"2"^(.1){\,}="3"
\ar@{}"1";"2"^(.8){\,}="4"
\ar@{=>}"3";"4"_{m}
&& x
}
\end{equation}
  where $m^{-1}$ is the horizontal inverse of the 2-arrow $m$. We
  denote by $\widetilde{m}: \al_\ast \Longrightarrow \be_\ast$
 the natural transformation which $m$ defines in this way. It is
 therefore an
 arrow \[ \widetilde{m}:\alpha_{\ast} \la \be_\ast\]
in the monoidal category $\mathcal{E}q(\gc)_{U}$.
With
  this notation, it follows that equation \eqref{defnualg} conjugates
  according to the scheme
\[\xymatrix@=70pt{
*+[F]{\gc_{U_{ijkl}}} 
\ar@/^2pc/[rr]^{i_{\lambda_{ij}(g_{jkl})}\,\,i_{ g_{ijl}}\,\, \lambda_{il} } 
 \ar@/_2pc/[rr]_{\lambda_{ij}\,\,
\lambda_{jk}\,\,
\lambda_{kl}} 
\ar@/^2pc/[rr]^(.3){}="1" \ar@/^2pc/[rr]^(.7){}="2"
  \ar@/_2pc/[rr]^(.3){}="3"\ar@/_2pc/[rr]^(.7){}="4"
 && *+[F]{\gc_{U_{ijkl}}}
\ar@{}"1";"3"^(.2){\,}="5"
\ar@{}"1";"3"^(.8){\,}="6"
\ar@{=>}"5";"6"_{3}
\ar@{}"2";"4"^(.2){\,}="7"
\ar@{}"2";"4"^(.8){\,}="8"
\ar@{=>}"7";"8"^{4}
\ar@{}"5";"6"^(.5){\,}="9"
\ar@{}"7";"8"_(.5){\,}="10"
\ar@{}"9";"10"^(.3){\,}="11"
\ar@{}"9";"10"^(.7){\,}="12"
\ar@{=}"11";"12"
}\] 
to the following  equation between the arrows  ``3'' and ``4''
 in the category $\mathcal{E}q(\gc)_{U_{ijkl}}$:
\begin{equation}
\label{eq:inu}
\widetilde{m}_{ijk}\,\,
{}^{g_{ijk}\,}\!\widetilde{m}_{ikl}\,\, j(\nu_{ijkl})
= (\lambda_{ij}\,\,\widetilde{m}_{jkl})\,\,
{}^{\lambda_{ij}(g_{jkl})\,\,}\!
\widetilde{m}_{ijl}
\end{equation} 
In such an equation, the term $j(\nu_{ijkl})$ is  the image of the
element  $\nu_{ijkl} \in \mathrm{Ar}(\gc)$ under
 inner conjugation functor\footnote{which
should not be confused with the inner conjugation homomorphism $i: G
\la \mathrm{Aut}(G)$
\eqref{def:inconj} which arises, as we saw at the end of \S 5.2, when $\gc$ is the stack
$\mathrm{Bitors}(G)$ associated to a bundle of groups $G$.
} 
\begin{equation}
\label{innerg}
 \gc \stackrel{j}{\la} \mathcal{E}q(\gc) 
\end{equation} 
 associated to  the group-like 
monoidal stack $\gc$. By an expression such as
${}^{g_{ijk\,}}\!\widetilde{m}_{ikl}$, we mean the conjugate
 of the 1-arrow   $\widetilde{m}_{ikl}$
 by  the
object  $j(g_{ijk})$
  in the monoidal category
$\mathcal{E}q(\gc)_{U_{ijkl}}$ . We observe here that the right
whiskerings of a 2-arrow $m$ or $\nu$ ({\it i.e.} the composition a
2-arrow with a 1-arrow which precedes it) have no significant effect
upon the conjugation operation which associates to a 2-arrow $m$
$(resp.\ \nu)$ in
$\pc$ the corresponding  natural transformation $\widetilde{m}$
 $(resp.\ j(\nu))$, an arrow in $\mathcal{E}q(\gc)$.
   This is why it was  harmless to ignore the right whiskerings in
 formula \eqref{defnualg},  and we will do so in similar contexts in the
 sequel.

\bigskip

Let us   display once more  the  cube \eqref{imabianchicube}, but now
decorated according to the  same
conventions  as in \eqref{imacube-dot}:
 \begin{equation}
 \label{imabianchicube1}
         \xymatrix@=15pt{
    &&x_i \ar[rrrrr]^{g_{ijl}} 
   \ar@2{.>} '[dd][ddddd]^(.3){g_{ikl}}
    \ar@{}[ddddd]_(.8){\,}="30"
    \ar@{}[rrrrr]_(.2){\,}="1"
    \ar@{}[ddddd]^(.3){\,}="2"
     &&&&& x_i
    \ar[ddddd]^{\lambda_{ij}(g_{jkl})}\\
    &&&&&&&\\
  *+[F]{x_i} \ar@{-->}[rrrrr]^(.65){g_{ijm}}
    \ar@{}[rrrrr]^(.3){\,}="25"
     \ar@2{.>}[uurr]^{g_{ilm}}
     \ar@{}[uurr]^(.3){\,}="26"
     \ar@{}[rrrrr]^(.7){\,}="5"
    \ar[ddddd]_{g_{ikm}}
    \ar@{}[ddddd]^(.6){\,}="9"
    &&&&&x_i
    \ar[uurr]^{\lambda_{ij}(g_{jlm})}
     \ar@{}[uurr]_(.3){\,}="22"
      \ar@{-->}[ddddd]^(.4){\lambda_{ij}(g_{jkm})}
    \ar@{}[ddddd]^(.3){\,}="21"
    \ar@{}[ddddd]_(.8){\,}="17"
    &&\\
    &&&&&&&\\
    &&&&&&&\\
    && x_i
    \ar@2{.>} '[rrr]^(.7){g_{ijk}} [rrrrr]
    \ar@{}[rrrrr]_(.3){\,}="14"
    &&
    &&&
  *+[F]{x_i}
    \\
    &&&&&&&
    \\ x_i
    \ar[uurr]_(.8){\lambda_{ik}(g_{klm})}
     \ar@{}[uurr]_(.4){\,}="34"
      \ar@{}[uurr]^(.7){\,}="29"
    \ar[rrrrr]_{g_{ijk}}
     \ar@{}[rrrrr]^(.3){\,}="33"
      \ar@{}[rrrrr]^(.8){\,}="18"
     &&&&&x_i
      \ar@{-->}[uurr]_(.6){\lambda_{ij}\lambda_{jk}(g_{klm})}
      && 
    \ar@{}"1";"2"^(.3){\,}="3"
     \ar@{}"1";"2"^(.8){\,}="4"
   \ar@{=>}"3";"4"^{\nu_{ijkl}}
   \ar@{}"25";"26"^(.2){\,}="27"
   \ar@{}"25";"26"^(.8){\,}="28"
   \ar@{=>}"27";"28"_{\nu_{ijlm}}
   \ar@{}"29";"30"^(.3){\,}="31"
   \ar@{}"29";"30"^(.7){\,}="32"
    \ar@{=>}"31";"32"^{\nu_{iklm}}
    \ar@{}"21";"22"^(.3){\,}="23"
    \ar@{}"21";"22"^(.7){\,}="24"
    \ar@{=>}"23";"24"_{\lambda_{ij}(\nu_{jklm})}
     \ar@{}"17";"18"^(.2){\,}="19"
       \ar@{}"17";"18"^(.8){\,}="20"
        \ar@{=>}"19";"20"_{\nu_{ijkm}}
         \ar@{}"33";"34"^(.05){\,}="35"
         \ar@{}"33";"34"^(.95){\,}="36"
     \ar@{=>}"35";"36"_{\{\widetilde{m}_{ijk},\,g_{klm}\}^{-1}}
        } 
\end{equation}  
 The commutativity of this diagram of 2-arrows  translates (according to the
 recipe  which  produced the algebraic equation \eqref{defnualg} from the
 cube
\eqref{imacube-dot})
 to the following 
 very twisted 
  3-cocycle condition for $\nu$
\footnote{\label{nu} 
This is essentially the
   3-cocycle equation (4.2.17) of \cite{2-gerbe}, but  with the terms in
   opposite order due to the fact that the somewhat imprecise
   definition of a 2-arrow $\nu$  given on page 71  of \cite{2-gerbe}  
yields  the inverse of the  2-arrow $\nu$ defined here by
   equation \eqref{defnualg}.
}: 
\begin{equation}
\label{3coceq}
\nu_{ijkl} \, 
({}^{\lambda_{ij}(g_{jkl})\,}\!\nu_{ijlm})\,\, \lambda_{ij}(\nu_{jklm}) \,\, 
 =  {}^{g_{ijk}\,}\!\nu_{iklm} \,\,
\{\widetilde{m}_{ijk},g_{klm}\}^{-1}\,\,
 (^{\lambda_{ij}\lambda_{jk}(g_{klm})\,}\!\nu_{ijkm})
\end{equation}

\noindent This is an equation satisfied
 by elements with values in
 $\mathrm{Ar}\,(\gc_{U_{ijklm}})$.
Note the  occurrence here of the term 
$\{\widetilde{m}_{ijk},g_{klm}\}^{-1}$,
corresponding to the lower face of \eqref{imabianchicube1}. While such
a term  does not exist in the standard definition of an 
 abelian \v{C}ech 3-cocycle equation,  
 non-abelian 3-cocycle relations
of this type  go back to  the
work of   P.  Dedecker \cite{ded}. They arise there   in  the context
of group 
  rather than   \v{C}ech cohomology, with his cocycles taking their
  values 
in an unnecessarily 
 restrictive precursor
  of a crossed square, which he calls a super-crossed group.

\bigskip

The following definition, which summarizes the previous discussion,
may be also viewed as a categorification of the notion of a $G$-valued
cocycle pair, as defined by equations \eqref{defgerbecoc}:

\begin{definition}
Let $\gc$ be a group-like monoidal stack on a space $X$, and
$\mathcal{U}$ an open covering of $X$. A $\gc$-valued  \v{C}ech
1-cocycle
 quadruple
is a quadruple of elements
 \begin{equation}
\label{alg3coc}
 (\lambda_{ij},\, \widetilde{m}_{ijk},\,g_{ijk},\, \nu_{ijkl})
\end{equation}
 satisfying the following conditions.
 The term  $\lambda_{ij} $ is an object in  the monoidal
 category  $\mathcal{E}q_{U_{ij}}(\gc_{|U_{ij}})$  and 
 $\widetilde{m}_{ijk}$ is  an arrow 
\begin{equation}
\label{def:mtilde1}
\widetilde{m}_{ijk}: j(g_{ijk})\,\lambda_{ik} \la \lambda_{ij} \, 
\lambda_{jk}
\end{equation}
in the corresponding monoidal category
$\mathcal{E}q_{U_{ijk}}(\gc_{|U_{ijk}})$.
Similarly,  $g_{ijk}$ is
 an object in the monoidal category $\gc_{U_{ijk}}$ and  
\[\nu_{ijkl}:\lambda_{ij}(g_{jkl}) \, g_{ijl} \la g_{ijk}\, g_{ikl}\]
an arrow  \eqref{def:nuijkl} in the corresponding monoidal category
$\gc_{U_{ijkl}}$. Finally we require that the two equations
 \eqref{eq:inu} and \eqref{3coceq}, which we reproduce here for the
 reader's convenience, be satisfied:
\begin{equation}
\label{cocsum}
\left \{
\begin{array}{lcl}
\widetilde{m}_{ijl}\,\,
{}^{g_{ijk}\,}\!\widetilde{m}_{ikl}\,\, j(\nu_{ijkl})
&=& (\lambda_{ij}\,\,\widetilde{m}_{jkl})\,\,
{}^{\lambda_{ij}(g_{jkl})\,\,}\!
\widetilde{m}_{ijl}
\\
\nu_{ijkl} \, 
({}^{\lambda_{ij}(g_{jkl})\,}\!\nu_{ijlm})\,\, \lambda_{ij}(\nu_{jklm}) \,\, 
& =&  {}^{g_{ijk}\,}\!\nu_{iklm} \,\,
\{\widetilde{m}_{ijk},g_{klm}\}^{-1}\,\,
 (^{\lambda_{ij}\lambda_{jk}(g_{klm})\,}\!\nu_{ijkm})
\end{array}
\right.
\end{equation}
\end{definition}

\bigskip

 Returning to our discussion, let us consider such a  $\gc$-valued \v{C}ech
 1-cocycle quadruple  \eqref{alg3coc}.
In order to transform the categorical  crossed module 
\eqref{innerg} into a weak analogue of a  
 crossed square, it is expedient for us  to restrict
ourselves, in both  the categories $\gc $ and $\mathcal{E}q(\gc)$, 
  to those arrows whose
source is the identity object. Diagram \eqref{innerg} then becomes
\begin{equation}
\label{crsq}
\xymatrix{
\mathrm{Ar}_I\,\gc \ar[d]_t \ar[r]^(.4)j & \mathrm{Ar}_I\,\mathcal{E}q(\gc)
\ar[d]^{t}\\
\mathrm{Ob} \, \gc \ar[r]_(.4)j & \mathrm{Ob}\, \mathcal{E}q(\gc) 
}
\end{equation}
where $t$ is the target map and the same symbol $j$ describes  the
components on objects and on arrows of the inner conjugation functor \eqref{innerg}.
 Recall  that one can assign to  any 
 arrow 
$u:X \la Y$  in a
 group-like  monoidal category
 the arrow $uX^{-1}:I \la YX^{-1}$ sourced at the identity, without loosing any
significant information. In particular, the arrow
$\widetilde{m}_{ijk}$ \eqref{def:mtilde1}
may  be  replaced by an arrow 
\[I \la \lam_{ij}\, \lam_{jk}\, \lam_{ik}^{-1}
\, j(g_{ijk})^{-1}\]
in $(\mathrm{Ar}_I\,\mathcal{E}q(\gc))_{U_{ijk}}$
 and the arrow $\nu_{ijkl}$  \eqref{def:nuijkl}
 by an arrow \[I \la 
g_{ijk}\,\,g_{ikl}\,\,g_{ijl}^{-1}\,\,\lambda_{ij}(g_{ijk})^{-1}\,,\]
 in $(\mathrm{Ar}_I\,\gc)_{U_{ijkl}}$
  which we again respectively denote by $\widetilde{m}_{ijk}$ and
$\nu_{ijkl}$. 
Our  quadruple \eqref{alg3coc}
then takes its values  in the weak square 
\begin{equation}
\label{crsq1}
\xymatrix{
(\mathrm{Ar}_I\,\gc)_{U_{ijkl}} \ar[d]_t \ar[r]^(.4)j 
& (\mathrm{Ar}_I\,\mathcal{E}q(\gc))_{U_{ijk}}
\ar[d]^{t}\\
(\mathrm{Ob} \, \gc)_{U_{ijk}} \ar[r]_(.4)j & 
(\mathrm{Ob}\, \mathcal{E}q(\gc))_{U_{ij}} 
}
\end{equation}
in the positions 
\begin{equation}
\label{crsq2}
\left(
\begin{matrix}
 \nu_{ijkl} & \widetilde{m}_{ijk}\\
g_{ijk} & \lambda_{ij}
\end{matrix}
\right) 
\end{equation}
Since  the evaluation action of $\mathcal{E}q(\gc)$ on $\gc$ produces a map 
\[\mathrm{Ar}_I\,\mathcal{E}q(\gc) \times \mathrm{Ob} \, \gc 
\la \mathrm{Ar}_I\,\gc
\]
which is the analog of the morphism \eqref{crhommap},
the quadruple  \eqref{alg3coc} may now  be viewed as 
a cocycle with values in what might be termed the (total complex associated to the)
  weak crossed square \eqref{crsq}. We will say that this modified quadruple
\eqref{crsq2} is a \v{C}ech 1-cocycle\footnote{ This
was called a 3-cocycle in
 \cite{2-gerbe}, but the present terminology is more appropriate.} for the covering
$\mathcal{U}$ on $X$
with values in   the (weak) crossed square \eqref{crsq}. Because of
the position of the different terms of the quadruple \eqref{crsq2} 
in the square \eqref{crsq1}, this terminology is consistent with
the fact that the  component $\nu_{ijkl}$ of such a 1-cocycle 
\eqref{alg3coc}  satisfies a sort of 3-cocycle
relation \eqref{3coceq}. 
The discussion in paragraph 6.2 will now be summarized  as follows in purely
algebraic terms:
\begin{proposition}
\label{prop:coc}
To a $\gc$-2-gerbe $\pc$ on $X$, locally trivialized by the choice of
objects $x_i$ in $\pc_{U_i}$ and local paths $\phi_{ij}$ \eqref{2phij},
is associated a 1-cocycle \eqref{alg3coc} with values in the weak
crossed square \eqref{crsq}.
\end{proposition}

\begin{remark}
{
\rm  
When $\gc$ is the $gr$-stack associated to a crossed module $\delta:
G \la \pi$,
this coefficient crossed module of $gr$-stacks is a stackified version
of the following crossed square associated by K.J.Norrie (see
\cite{nor}, \cite{brown} theorem 3.5)
 to the  crossed
module $G {\la} \pi$:

\begin{equation}
\label{diagcrmod}
  \xymatrix@=12pt{
G \ar[d]_{\delta}  \ar[rr] && \mathrm{Der}^{\ast}(\pi, G) \ar[d]\\
\pi \ar[rr] && \mathrm{Aut}(G \rightarrow   \pi)
}\end{equation}

\medskip

\noindent It is however less restrictive than Norrie's version,
 since the latter
  corresponds to the  diagram of $gr$-stacks 
\[ \gc \la \mathrm{Isom} (\gc)
\]
\noindent whereas we really need to consider, as in 
\eqref{innerg}, self-equivalences of the monoidal stack $\gc$,
rather than automorphisms. To phrase it differently, we need to
replace  the term
$\mathrm{Aut}(G \la \pi)$ in the square \eqref{diagcrmod} by the weak
automorphisms of the crossed module $G \la \pi$,
 as discussed in remark \ref{rem:crmod},
and modify the set of crossed homomorphisms
 $\mathrm{Der}^{\ast}(\pi, G)$ accordingly.
}
\end{remark}

\subsection    \noindent {\bf Coboundary relations}

 \bigskip

 We now choose a second set of local objects $x'_{i} \in
\pc_{U_i}$, and  of local arrows \eqref{2phij}
\[ \phi'_{ij}: x'_j \la x'_i\,\] 
By proposition \ref{prop:coc},
 these determine a second  crossed square valued 1-cocycle
 \begin{equation}
\label{alg3coc1}
(\lambda'_{ij},\, \widetilde{m}'_{ijk},\, g'_{ijk},\, \nu'_{ijkl})\,.
\end{equation}
In order to compare it with the 1-cocycle  \eqref{alg3coc}, we proceed
as we did   in section 5.2 above, but now in a 2-categorical setting. 
We choose once more an arrow $\chi_i$ \eqref{def:xi}.
There now  exist  1-arrows  $\te_{ij}$, and  2-arrows $\zeta_{ij}$
in $\pc_{U_{ij}}$.
\begin{equation}
\label{def:deij2}
\xymatrix@R=15pt@C=40pt{
x_j \ar[r]^{\phi_{ij}} \ar[dd]_{\chi_j} & x_i \ar[d]^{\chi_i} \\
& x'_i  \ar@{}[d]_(.01){\,}="1" \ar[d]^{\te_{ij}} \\
x'_j \ar[r]_{\phi'_{ij}}
\ar@{}[r]^{\,}="2"
 & x'_i
\ar@{}"1";"2"_(.2){\,}="3"
\ar@{}"1";"2"_(.8){\,}="4"
\ar@{=>}"4";"3"^{\zeta_{ij}}
}
\end{equation} 
 The arrow $\chi_i$ induces by 
conjugation a self-equivalence $r_i: \gc \la \gc$ and  2-arrows 
\begin{equation}
\label{def-mui1}
\xymatrix@R=30pt@C=40pt{ x_i \ar[d]_{\chi_i} 
\ar[r]^u & x_i   \ar[d]^{\chi_i}  \ar@{}[d]_{\,}="1"
\\
x'_i  \ar[r]_{r_i(u)} \ar@{}[r]_{\,}="2"
 & x'_i
\ar@{}"1";"2"_(.2){\,}="3"
\ar@{}"1";"2"_(.8){\,}="4"
\ar@{=>}"3";"4"
}
\end{equation}
which  are functorial in $u$.
 Furthermore, the  diagram \eqref{def:deij2} induces by conjugation a
diagram in $\gc_{U_{ij}}$:
\begin{equation}
\label{def:deij3}
\xymatrix@R=15pt@C=40pt{
\gc \ar[r]^{\lambda_{ij}} \ar[dd]_{r_j} & \gc \ar[d]^{r_i} \\
& \gc  \ar@{}[d]_(.01){\,}="1" \ar[d]^{i_{\te_{ij}}} \\
\gc \ar[r]_{\lambda'_{ij}}
\ar@{}[r]^{\,}="2"
 & \gc
\ar@{}"1";"2"_(.2){\,}="3"
\ar@{}"1";"2"_(.8){\,}="4"
\ar@{=>}"4";"3"^{\widetilde{\zeta}_{ij}}
}
\end{equation} 
with
$\widetilde{\zeta}_{ij}$
the natural transformation induced by $\zeta_{ij}$.
Consider now the diagram of 2-arrows
\begin{equation}
\label{diagcoboun1a}
\xymatrix@R=15pt@C=50pt{
&x_k \ar[ddl]_{\phi_{jk}} \ar[rr]^{\phi_{ik}}
\ar@{}[rr]_(.6){\,}="1"
\ar@{}[rr]_(.2){\,}="25"
\ar '[dd]_(.5){\chi_{k}} '[dddd]   [ddddd]
\ar@{}[dddd]^(.45){\,}="17"
\ar@{}[dddd]^(.25){\,}="26"
   &&x_i \ar[d]^{\chi_i}
\ar@{}[d]_(.5){\,}="9"
\ar@{}[ddl]_(.99){\,}="10"
  \ar[ddl]_{g_{ijk}}\\
&&&  x'_i \ar[ddl]^(.7){r_i(g_{ijk}) } \ar[dddd]^{\te_{ik}} \\
   x_j \ar[rr]_(.65){\phi_{ij}}
\ar@{}[rr]^(.6){\,}="2"
\ar[dd]_{\chi_{j}}
\ar@{}[dd]^(.9){\,}="18"
&&  x_i
\ar[d]_(.45){\chi_i}
& \\
&& x'_i \ar[d]^{\te_{ij}} \ar@{}[d]_(.01){\,}="13"
&\\
 x'_j  
\ar [rr]^(.65){\phi'_{ij}}\ar@{}[rr]^(.65){\,}="14"
\ar@{}[rr]_(.35){\,}="21"
 \ar[ddd]_{\te_{jk}}
\ar@{}[ddd]^(.25){\,}="22"
 & & x'_i
\ar[ddd]^(.25){ \lambda'_{ij}(\te_{jk})}
 & \\
& x'_k  \ar[ddl]_{\phi'_{jk}} 
\ar '[r]^(.5){\phi'_{ik}} [rr]
\ar@{}[rr]_(.6){\,}="5"
 && x'_i \ar[ddl]^{g'_{ijk}} \\
&& &\\
  x'_j \ar[rr]_ {\phi'_{ij}}
\ar@{}[rr]_(.6){\,}="6"
 &&  x'_i & 
\ar@{}"1";"2"^(.4){\,}="3"
 \ar@{}"1";"2"^(.7){\,}="4"
 \ar@{=>}"3";"4"_{m_{ijk}}
\ar@{}"5";"6"^(.4){\,}="7"
 \ar@{}"5";"6"^(.7){\,}="8"
 \ar@{=>}"7";"8"_{m'_{ijk}}
\ar@{}"9";"10"^(.4){\,}="11"
\ar@{}"9";"10"^(.7){\,}="12"
 \ar@{=>}"11";"12"
 \ar@{}"13";"14"^(.2){\,}="15"
 \ar@{}"13";"14"^(.7){\,}="16"
 \ar@{=>}"16";"15"^{\zeta_{ij}}
 \ar@{}"17";"18"^(.4){\,}="19"
\ar@{}"17";"18"^(.7){\,}="20"
 \ar@{=>}"19";"20"_{\zeta_{jk}}
\ar@{}"21";"22"^(.3){\,}="23"
\ar@{}"21";"22"^(.8){\,}="24"
\ar@{=>}"23";"24"
\ar@{}"25";"26"^(.2){\,}="27"
\ar@{}"25";"26"^(.8){\,}="28"
\ar@{=>}"28";"27"^{\zeta_{ik}}
}
\end{equation}  
which extends \eqref{diagcoboun}.
Three of its 2-arrows are of the form $\zeta_{ij}$, the top and the
bottom ones  are of the form $m_{ijk}$ \eqref{def:mijk}. The unlabelled lower
front 2-arrow and the right-hand upper are respectively part of the definitions 
of $\lambda'_{ij}(\te_{jk})$ and of $r_i(g_{ijk})$. Since  these seven
2-arrows are invertible, diagram
\eqref{diagcoboun1a}  uniquely defines a 2-arrow $b_{ijk}$ filling in the
remaining lower right-hand square:
\begin{equation}
\label{def:bijk}
\xymatrix@C=40pt{
x'_i 
\ar@{}[d]_(.5){\,}="2"
\ar@{}
\ar[d]_{\te_{ij}}
 &  \ar[l]_{r_i(g_{ijk})} x'_i \ar[dd]^{\te_{ik}}
\ar@{}[d]_(.9){\,}="1"
\\
x'_i
 \ar[d]_{\lambda'_{ij}(\te_{jk})} &\\
x'_i &x'_i \ar[l]^{g'_{ijk}}
\ar@{}"1";"2"^(.3){\,}="3"
 \ar@{}"1";"2"^(.6){\,}="4"
 \ar@{=>}"3";"4"_{b_{ijk}}
}
\end{equation}
so that diagram \eqref{diagcoboun1a} becomes the following  commutative diagram
of 2-arrows, which we directly display  in decorated form, according
to the conventions of \eqref{imacube-dot}: 
\begin{equation}
\label{diagcoboun2}
\xymatrix@R=15pt@C=50pt{
&*+[F]{x_k} \ar@2{.>}[ddl]_{\phi_{jk}} \ar[rr]^{\phi_{ik}}
\ar@{}[rr]_(.6){\,}="1"
\ar@{}[rr]_(.2){\,}="25"
\ar@{-->}[ddddd]_(.25){\chi_{k}}
\ar@{}[dddd]^(.45){\,}="17"
\ar@{}[dddd]^(.25){\,}="26"
   &&x_i \ar[d]^{\chi_i}
\ar@{}[d]_(.5){\,}="9"
\ar@{}[ddl]_(.99){\,}="10"
  \ar[ddl]_{g_{ijk}}\\
&&&  x'_i \ar[ddl]^(.7){r_i(g_{ijk}) } \ar[dddd]^{\te_{ik}}
\ar@{}[dddd]_(.8){\,}="31"
 \\
   x_j \ar@2{.>}[rr]_(.65){\phi_{ij}}
\ar@{}[rr]^(.6){\,}="2"
\ar[dd]_{\chi_{j}}
\ar@{}[dd]^(.9){\,}="18"
&&  x_i
\ar@2{.>}[d]_(.45){\chi_i}
& \\
&& x'_i \ar@2{.>}[d]^{\te_{ij}} \ar@{}[d]_(.01){\,}="13"
\ar@{}[d]_(.6){\,}="32"
&\\
 x'_j  
\ar[rr]^(.65){\phi'_{ij}}\ar@{}[rr]^(.65){\,}="14"
\ar@{}[rr]_(.35){\,}="21"
 \ar[ddd]_{\te_{jk}}
\ar@{}[ddd]^(.25){\,}="22"
 & & x'_i
\ar@2{.>}[ddd]^(.25){ \lambda'_{ij}(\te_{jk})}
 & \\
& x'_k  \ar[ddl]_{\phi'_{jk}} 
\ar@{-->} '[r]^(.5){\phi'_{ik}} [rr]
\ar@{}[rr]_(.6){\,}="5"
 && x'_i \ar@{-->}[ddl]^{g'_{ijk}} \\
&& &\\
  x'_j \ar[rr]_ {\phi'_{ij}}
\ar@{}[rr]_(.6){\,}="6"
 &&  *+[F]{x'_i} & 
\ar@{}"1";"2"^(.4){\,}="3"
 \ar@{}"1";"2"^(.7){\,}="4"
 \ar@{=>}"3";"4"_{m_{ijk}}
\ar@{}"5";"6"^(.4){\,}="7"
 \ar@{}"5";"6"^(.7){\,}="8"
 \ar@{=>}"7";"8"_{m'_{ijk}}
\ar@{}"9";"10"^(.4){\,}="11"
\ar@{}"9";"10"^(.7){\,}="12"
 \ar@{=>}"11";"12"
 \ar@{}"13";"14"^(.2){\,}="15"
 \ar@{}"13";"14"^(.7){\,}="16"
 \ar@{=>}"16";"15"^{\zeta_{ij}}
 \ar@{}"17";"18"^(.4){\,}="19"
\ar@{}"17";"18"^(.7){\,}="20"
 \ar@{=>}"19";"20"_{\zeta_{jk}}
\ar@{}"21";"22"^(.3){\,}="23"
\ar@{}"21";"22"^(.8){\,}="24"
\ar@{=>}"23";"24"
\ar@{}"25";"26"^(.2){\,}="27"
\ar@{}"25";"26"^(.8){\,}="28"
\ar@{=>}"28";"27"^{\zeta_{ik}}
\ar@{}"31";"32"^(.3){\,}="33"
 \ar@{}"31";"32"^(.6){\,}="34"
 \ar@{=>}"33";"34"_{b_{ijk}}
}   
\end{equation}
We derive from this
diagram the algebraic  equation
\[
(\lambda'_{ij}(\te_{kl}) \ast \zeta_{ij}) \,\, (\phi'_{ij} \ast
\zeta_{jk})\,\,m'_{ijk} =((\lam'_{ij}(\te_{jk})\, \te_{ij}\,\chi_i)\ast
m_{ijk})\,\,
\,b_{ijk}\, \,(g'_{ijk}\ast \zeta_{ik})
\]
for the equality between the two corresponding 2-arrows between the
decorated paths.
With the same notations as for equation \eqref{eq:inu}, the conjugated
version of this   equation is
\begin{equation}
\label{cobouneq1}
{}^{\lambda'_{ij}(\te_{jk})\,}\!\widetilde{\zeta}_{ij} \,\, 
\lambda'_{ij}(\widetilde{\zeta}_{jk})\, \widetilde{m}'_{ijk}
=
{}^{\lambda'_{ij}(\te_{jk})\,\te_{ij}\,r_i} \,\widetilde{m}_{ijk}\,\,
j(b_{ijk}) \,\, ({}^{g'_{ijk}\,}\!\widetilde{\zeta}_{ik})
\end{equation} 
It is 
the analogue, with the present conventions, of  equation
 \cite{2-gerbe} (4.4.12).

\bigskip 

A  second coboundary condition relates the cocycle  quadruples 
\eqref{alg3coc} and \eqref{alg3coc1}. In geometric terms, it asserts
the commutativity of the following diagram of 2-arrows, in which the
unlabelled 2-arrow in the middle of the right vertical face is
 $\{\widetilde{\zeta}_{ij},\, g_{jkl}\}^{-1}$ defined in the same way
 as the 2-arrow  which we denoted $\{\widetilde{m}_{ijk},\,g_{klm}\}$ 
\eqref{def:brack}:
\begin{equation}
\label{diagcoboun2long}
\xymatrix@R=15pt@C=70pt{
&x_i \ar[dl]_{g_{ikl}} \ar[rr]^{g_{ijl}} 
 \ar '[d]^(.5){\chi_i} [dd]
\ar@{}[rr]_(.6){\,}="41"
&& x_i \ar[dd]^{\chi_i}
\ar[dl]^{\lambda_{ij}(g_{jkl})}\\
x_i  \ar@{}[rr]^(.6){\,}="42"
 \ar[ddd]_(.4){\chi_i}
\ar[rr]_(.65){g_{ijk}}&&x_i  \ar[ddd]_{\chi_i}  &\\
&x'_i \ar[ddl]_{r_i(g_{ikl})} \ar '[r]^(.5){r_i(g_{ijl})} [rr]
\ar@{}[rr]_(.6){\,}="1"
\ar@{}[rr]_(.2){\,}="25"
\ar '[dd]_(.45){\te_{il}} '[dddd] [ddddd]
\ar@{}[dddd]^(.45){\,}="17"
\ar@{}[dddd]^(.25){\,}="26"
   &&x'_i
\ar[d]^{\te_{ij}}
\ar@{}[d]_(.5){\,}="9"
\ar@{}[ddl]_(.99){\,}="10"
  \ar[ddl]_(.4){r_i(\lambda_{ij}(g_{jkl}))\!\!}\\
&&&  x'_i \ar[ddl]^(.7){\lambda'_{ij}(r_j(g_{jkl})) }
\ar[dddd]^{\lambda'_{ij}(\te_{jl})}
\ar@{}[dddd]_(.8){\,}="31"
 \\
   x'_i \ar[rr]_(.65){r_i(g_{ijk})}
\ar@{}[rr]^(.6){\,}="2"
\ar[dd]_{\te_{ik}}
\ar@{}[dd]^(.9){\,}="18"
&&  x'_i
\ar[d]_{\te_{ij}}
& \\
&& x'_i
 \ar[d]^(.45){\lambda'_{ij}(\te_{jk})}
 \ar@{}[d]_(.01){\,}="13"
\ar@{}[d]_(.6){\,}="32"
&\\
 x'_i  \ar[ddd]_{ \lambda'_{ik}(\te_{kl})}
\ar[rr]^(.65){g'_{ijk}}\ar@{}[rr]^(.65){\,}="14"
\ar@{}[rr]_(.45){\,}="21"
 \ar[ddd]_{}
\ar@{}[ddd]^(.35){\,}="22"
 & & x'_i \ar[ddd]_(.7){\lambda'_{ij}\lambda'_{jk}(\te_{kl})}
 & \\
& x'_i  \ar[ddl]_{g'_{ikl}} 
 \ar '[r]_ {g'_{ijl}} [rr]
\ar@{}[rr]_(.3){\,}="5"
 && x'_i \ar[ddl]^{\lambda'_{ij}(g'_{jkl})} \\
&& &\\
  x'_i   \ar[rr]^(.3){g'_{ijk}}
\ar@{}[rr]_(.3){\,}="6"
 &&  x'_i & 
\ar@{}"1";"2"^(.4){\,}="3"
 \ar@{}"1";"2"^(.7){\,}="4"
 \ar@{=>}"3";"4"_{r_i(\nu_{ijkl})}
\ar@{}"5";"6"^(.4){\,}="7"
 \ar@{}"5";"6"^(.7){\,}="8"
 \ar@{=>}"7";"8"_{\nu'_{ijkl}}
\ar@{}"9";"10"^(.4){\,}="11"
\ar@{}"9";"10"^(.7){\,}="12"
 \ar@{=>}"11";"12"
 \ar@{}"13";"14"^(.2){\,}="15"
 \ar@{}"13";"14"^(.7){\,}="16"
 \ar@{=>}"16";"15"^{b_{ijk}}
 \ar@{}"17";"18"^(.4){\,}="19"
\ar@{}"17";"18"^(.7){\,}="20"
 \ar@{=>}"19";"20"_{b_{ikl}}
\ar@{}"21";"22"^(.38){\,}="23"
\ar@{}"21";"22"^(.68){\,}="24"
\ar@{=>}"24";"23"_{\:\:\;\{\widetilde{m}'_{ijk},\,\te_{kl}\}}
\ar@{}"25";"26"^(.2){\,}="27"
\ar@{}"25";"26"^(.8){\,}="28"
\ar@{=>}"28";"27"^{b_{ijl}}
\ar@{}"31";"32"^(.3){\,}="33"
 \ar@{}"31";"32"^(.6){\,}="34"
 \ar@{=>}"33";"34"_{\lap_{ij}(b_{jkl})}
\ar@{}"41";"42"^(.4){\,}="43"
 \ar@{}"41";"42"^(.7){\,}="44"
 \ar@{=>}"43";"44"_{\nu_{ijkl}}
}
\end{equation}
This  cubic diagram  
 compares the 2-arrows  $\nu_{ijkl}$ and $\nu'_{ijkl}$, which
are respectively its
 top and bottom faces. It actually consists of two separate 
 cubes. The upper one is
trivially commutative, as it simply defines the 2-arrow
$r_i(\nu_{ijkl})$, which is  the common face between the two
cubes considered. 

\begin{lemma}
\label{cobounlem}
The cube of 2-arrows \eqref{diagcoboun2long} is
commutative\footnote{ In   \cite{2-gerbe} \S 4.9  the
 corresponding assertion is implicitly
assumed to be true, although no proof is given there.}.
\end{lemma}

{\bf Proof:}  The proof that the full diagram
\eqref{diagcoboun2long} commutes is very similar to
the proof of  Proposition \ref{2com}. We consider a hypercube analogous
to diagram \eqref{imahypercube}, and which  therefore consists of eight
cubes called left, right, top, bottom, front, back, inner and
outer. The outer cube in this diagram is the cube
\eqref{diagcoboun2long}.
 We will  now
describe the seven other cubes. Since   these seven are
commutative,  this will suffice in order to prove that the outer one
also is, so that  the lemma will be proved. As this hypercubic diagram
is somewhat more complicated than \eqref{imahypercube}, we will 
describe it in words, instead of  displaying it.

\bigskip

\noindent The top cube is a copy of cube \eqref{imacube},
oriented so that its face $\nu_{ijkl}$ is on top,
consistently with the top face of \eqref{diagcoboun2long}. The  bottom
cube is a  cube of similar type  which defines
the 2-arrow $\nu'_{ijkl}$. Since it is built  from objects
$x'$, arrows 
$\phi'$ and $g'$ and 2-arrows $m'$ and $\nu'$, we will refer to it as
the primed version of \eqref{imacube} . It is 
 oriented so that   $\nu'_{ijkl}$ is the bottom face.

\bigskip

 \noindent We now
describe  the six other cubes. Four of these are of the type
 \eqref{diagcoboun2}. If we denote the latter by the symbol $P_{ijk}$
 determined by its indices, these are respectively  the left cube  
 $P_{ikl}$, the back cube  $P_{ijl}$,  the inner cube $P_{jkl}$
 and the front cube $P_{ijk}$.  Each of the first three rests on the
 corresponding face $m'_{ikl}$, $m'_{ijl}$, and  $m'_{jkl}$ of the bottom
 cube, and is attached at the top to the similar face $m$ of the top
 cube. 
The cube $P_{ijk}$  is  attached to 
 the corresponding  face $m_{ijk}$ of the top cube, but it does
 not constitute the full front cube. Below it is the following primed
 version of the cube \eqref{imabianchicube-bracket}, but
 now  associated to the face
 $\{\widetilde{m}'_{ijk}
,\, \te_{kl}\}$
 (rather
 than as in  \eqref{imabianchicube-bracket}  to the  face
 $\{\widetilde{m}_{ijk},\,g_{klm}\}$).
\begin{equation} 
 \label{imabianchicube-bracket2}
   \xymatrix@=15pt{
    &&x'_k\ar[rrrrr]^{\phi'_{jk}} 
   \ar '[dd][ddddd]^(.3){\phi'_{ik}}
 \ar@{}[rrrrr]_(.3){\,}="45"
    \ar@{}[ddddd]^(.3){\,}="46"
     &&&&& x'_j
\ar[ddddd]^{\phi'_{ij}}
\\
    &&&&&&&\\
    x'_k\ar[rrrrr]^(.65){\phi'_{jk}}
 \ar@{}[rrrrr]_(.3){\,}="50"
    \ar@{}[ddddd]^(.3){\,}="49"
     \ar[uurr]^{\te_{kl}}
\ar[ddddd]_{\phi'_{ik}}
    &&&&&x'_j
 \ar[uurr]^{\lambda'_{jk}(\te_{kl})}
    \ar[ddddd]^(.4){\phi'_{ij}}
    &&\\
    &&&&&&&\\
    &&&&&&&\\
    &&x'_i
    \ar '[rrr]^(.7){g'_{ijk}} [rrrrr]
    \ar@{}[rrrrr]_(.3){\,}="14"
    &&
    &&&
    x'_i
    \\
    &&&&&&&
    \\ x'_i 
    \ar[uurr]_(.75){\lambda'_{ik}(\te_{kl})}
      \ar@{}[uurr]^(.3){\,}="42"
    \ar[rrrrr]_{g'_{ijk}}
     \ar@{}[rrrrr]^(.4){\,}="41"
     &&&&& x'_i
      \ar[uurr]_(.6){\lambda'_{ij}\lambda'_{jk}(\te_{kl})}
      &&
\ar@{}"41";"42"^(.2){\,}="43"
 \ar@{}"41";"42"^(.7){\,}="44"
 \ar@{=>}"44";"43"^{\;\;\;\;\;\{\widetilde{m}_{ijk},\,\te_{kl}\}}
\ar@{}"45";"46"^(.2){\,}="47"
 \ar@{}"45";"46"^(.7){\,}="48"
 \ar@{=>}"48";"47"_{m_{ijk}}
\ar@{}"49";"50"^(.2){\,}="51"
 \ar@{}"49";"50"^(.7){\,}="52"
 \ar@{=>}"51";"52"_{m_{ijk}}
      }
\end{equation}

\vspace{.75cm}

Finally, the right cube is itself constituted of two cubes. The lower
one  constructs the 2-arrow
$\lambda'_{ij}(b_{jkl})$, starting from the 2-arrow $b_{jkl}$
 \eqref{def:bijk}.
 The upper one is
another  commutative cube of  type  \eqref{imabianchicube-bracket},
 but this time
associated to the face $\{\widetilde{\zeta}_{ij},\, g_{jkl}\}$.   This is
  the following commutative cube
whose four  unlabelled
2-arrows are the obvious ones.

\bigskip

\begin{equation} 
 \label{imabianchicube-bracket1}
   \xymatrix@R=15pt@C=38pt{
    &&x_j\ar[rrrrr]^{\phi_{ij}} 
   \ar '[dd][ddddd]^(.3){g_{jkl}}
    \ar@{}[ddddd]_(.8){\,}="30"
    \ar@{}[rrrrr]_(.2){\,}="1"
    \ar@{}[ddddd]^(.3){\,}="2"
     &&&&& x_i
    \ar[ddddd]^{\lambda_{ij}(g_{jkl})}
\\
    &&&&&&x'_i \ar@{<-}[ur]^{\chi_i}
\ar[ddddd]^{r_i(\lambda_{ij}(g_{jkl}))}
&\\
    x'_j\ar[rrrrr]^(.65){\phi'_{ij}}
    \ar@{}[rrrrr]^(.8){\,}="25"
     \ar@{<-}[uurr]^{\chi_{j}}
     \ar@{}[rrrrr]^(.7){\,}="5"
     \ar[ddddd]_{r_j(g_{jkl})}
    \ar@{}[ddddd]^(.6){\,}="9"
    &&&&&x_i
\ar@{}[uurr]^(.3){\,}="26"
 \ar@{<-}[ur]_{\te_{ij}}
\ar@{}[ddddd]^(.15){\,}="52"
     \ar@{}[uurr]_(.3){\,}="22"
    \ar[ddddd]_(.3){\lap_{ij}(r_j(g_{jkl}))}
\ar@{}[ddddd]_(.6){\,}="42"
    \ar@{}[ddddd]^(.3){\,}="21"
    \ar@{}[ddddd]_(.8){\,}="17"
    &&\\
    &&&&&&&\\
    &&&&&&&\\
    &&x_j
    \ar '[rrr]^(.7){\phi_{ij}}   '[rrrr] [rrrrr]
    \ar@{}[rrrrr]_(.3){\,}="14"
    &&
    &&&
    x_i
     \\
    &&&&&&x'_i  \ar@{<-}[ur]_(.6){\chi_{i}}
\ar[dl]^{\te_{ij}}
\ar@{}[uuuuu]^(.25){\,}="51"
 &
    \\ x'_j \ar@{}[rrrrr]_(.5){\,}="41"
    \ar@{<-}[uurr]_(.6){\chi_{j}}
     \ar@{}[uurr]_(.7){\,}="34"
      \ar@{}[uurr]^(.7){\,}="29"
    \ar[rrrrr]_{\phi'_{ij}}
     \ar@{}[rrrrr]^(.3){\,}="33"
      \ar@{}[rrrrr]^(.8){\,}="18"
     &&&&& x'_i
      && \Box
    \ar@{}"1";"2"^(.3){\,}="3"
     \ar@{}"1";"2"^(.8){\,}="4"
   \ar@{}"25";"26"^(.3){\,}="27"
   \ar@{}"25";"26"^(.7){\,}="28"
   \ar@{=>}"27";"28"^{\zeta_{ij}}
   \ar@{}"29";"30"^(.3){\,}="31"
   \ar@{}"29";"30"^(.7){\,}="32"
    \ar@{}"21";"22"^(.3){\,}="23"
    \ar@{}"21";"22"^(.7){\,}="24"
     \ar@{}"17";"18"^(.2){\,}="19"
       \ar@{}"17";"18"^(.8){\,}="20"
         \ar@{}"34";"14"^(.4){\,}="35"
         \ar@{}"34";"14"^(.7){\,}="36"
 \ar@{=>}"35";"36"_{\zeta_{ij}}
        \ar@{}"41";"42"^(.4){\,}="43"
\ar@{}"41";"42"^(.6){\,}="44"
\ar@{}"51";"52"^(.5){\,}="53"
\ar@{}"51";"52"^(.8){\,}="54"
\ar@{=>}"53";"54"_(.8){\{\widetilde{\zeta}_{ij},\,g_{jkl}\}}
 }
\end{equation}

\newpage

In order to translate the commutativity of the cube
\eqref{diagcoboun2long} into an algebraic expression, we now decorate it
as follows,
invoking once more the  conventions of  diagram \eqref{imacube-dot}:
\begin{equation}
\label{diagcoboun4long}
\xymatrix@R=15pt@C=70pt{
&*+[F]{x_i} \ar@2{.>}[dl]_{g_{ikl}} \ar[rr]^{g_{ijl}} 
 \ar@{-->} '[d]^(.5){\chi_i} [dd]
\ar@{}[rr]_(.6){\,}="41"
&& x_i \ar[dd]^{\chi_i}
\ar[dl]^{\lambda_{ij}(g_{jkl})}\\
x_i  \ar@{}[rr]^(.6){\,}="42"
 \ar[ddd]_(.4){\chi_i}
\ar@2{.>}[rr]_(.65){g_{ijk}}&&x_i  \ar@2{.>}[ddd]_{\chi_i}  &\\
&x'_i \ar[ddl]_{r_i(g_{ikl})} \ar '[r]^(.5){r_i(g_{ijl})} [rr]
\ar@{}[rr]_(.6){\,}="1"
\ar@{}[rr]_(.2){\,}="25"
\ar@{-->} '[dd]_(.45){\te_{il}} '[dddd] [ddddd]
\ar@{}[dddd]^(.45){\,}="17"
\ar@{}[dddd]^(.25){\,}="26"
   &&x'_i
\ar[d]^{\te_{ij}}
\ar@{}[d]_(.5){\,}="9"
\ar@{}[ddl]_(.99){\,}="10"
  \ar[ddl]_(.4){r_i(\lambda_{ij}(g_{jkl}))\!\!}\\
&&&  x'_i \ar[ddl]^(.7){\lambda'_{ij}(r_j(g_{jkl})) }
\ar[dddd]^{\lambda'_{ij}(\te_{jl})}
\ar@{}[dddd]_(.8){\,}="31"
 \\
   x'_i \ar[rr]_(.65){r_i(g_{ijk})}
\ar@{}[rr]^(.6){\,}="2"
\ar[dd]_{\te_{ik}}
\ar@{}[dd]^(.9){\,}="18"
&&  x'_i
\ar@2{.>}[d]_{\te_{ij}}
& \\
&& x'_i
 \ar@2{.>}[d]^(.45){\lambda'_{ij}(\te_{jk})}
 \ar@{}[d]_(.01){\,}="13"
\ar@{}[d]_(.6){\,}="32"
&\\
 x'_i  \ar[ddd]_{ \lambda'_{ik}(\te_{kl})}
\ar[rr]^(.65){g'_{ijk}}\ar@{}[rr]^(.65){\,}="14"
\ar@{}[rr]_(.45){\,}="21"
 \ar[ddd]_{}
\ar@{}[ddd]^(.35){\,}="22"
 & & x'_i \ar@2{.>}[ddd]_(.7){\lambda'_{ij}\lambda'_{jk}(\te_{kl})}
 & \\
& x'_i  \ar[ddl]_{g'_{ikl}} 
 \ar@{-->} '[r]_ {g'_{ijl}} [rr]
\ar@{}[rr]_(.3){\,}="5"
 && x'_i \ar@{-->}[ddl]^{\lambda'_{ij}(g'_{jkl})} \\
&& &\\
  x'_i   \ar[rr]^(.3){g'_{ijk}}
\ar@{}[rr]_(.3){\,}="6"
 && *+[F]{x'_i} & 
\ar@{}"1";"2"^(.4){\,}="3"
 \ar@{}"1";"2"^(.7){\,}="4"
 \ar@{=>}"3";"4"_{r_i(\nu_{ijkl})}
\ar@{}"5";"6"^(.4){\,}="7"
 \ar@{}"5";"6"^(.7){\,}="8"
 \ar@{=>}"7";"8"_{\nu'_{ijkl}}
\ar@{}"9";"10"^(.4){\,}="11"
\ar@{}"9";"10"^(.7){\,}="12"
 \ar@{=>}"11";"12"
 \ar@{}"13";"14"^(.2){\,}="15"
 \ar@{}"13";"14"^(.7){\,}="16"
 \ar@{=>}"16";"15"^{b_{ijk}}
 \ar@{}"17";"18"^(.4){\,}="19"
\ar@{}"17";"18"^(.7){\,}="20"
 \ar@{=>}"19";"20"_{b_{ikl}}
\ar@{}"21";"22"^(.38){\,}="23"
\ar@{}"21";"22"^(.68){\,}="24"
\ar@{=>}"24";"23"_{\:\:\;\{\widetilde{m}'_{ijk},\,\te_{kl}\}}
\ar@{}"25";"26"^(.2){\,}="27"
\ar@{}"25";"26"^(.8){\,}="28"
\ar@{=>}"28";"27"^{b_{ijl}}
\ar@{}"31";"32"^(.3){\,}="33"
 \ar@{}"31";"32"^(.6){\,}="34"
 \ar@{=>}"33";"34"_{\lap_{ij}(b_{jkl})}
\ar@{}"41";"42"^(.4){\,}="43"
 \ar@{}"41";"42"^(.7){\,}="44"
 \ar@{=>}"43";"44"_{\nu_{ijkl}}
}
\end{equation}
Reading  off the
  two  composite  2-arrows between the decorated 1-arrows 
\eqref{diagcoboun4long}, and taking  into account the appropriate
 whiskerings, we see that the commutativity of diagram
 \eqref{diagcoboun4long}
 is equivalent to the following algebraic equation.
\begin{equation}
\label{cobouneq2}
\begin{split}
({}^{\lap_{ij}\lap_{jk}(\te_{kl})\lap_{ij}(\te_{jk})\te_{ij}
r_i \;\;}\!\nu_{ijkl})\,
\,&({}^{\lap_{ij}\lap_{jk}(\te_{kl})\lap_{ij}(\te_{jk})
\,\,}\!\{\widetilde{\zeta}_{ij},\,g_{jkl}\}^{-1}) \,\,
\lap_{ij}(b_{jkl})\,\,({}^{\lap_{ij}(g'_{jkl})\,}\!b_{ijl})\\
&=
({}^{\lap_{ij}\lap_{jk}(\te_{kl})\,}\!b_{ijk})\,
\{\widetilde{m}'_{ijk},\,\te_{kl}\}
\,({}^{g'_{ijk}\,}\!b_{ikl})\, \,\nu'_{ijkl}
\end{split}
\end{equation}
This equation is the analogue, under  our present conventions, of equation
(4.4.15) of \cite{2-gerbe}. It
describes the manner in which the various  terms  of type $b_{ijk}$ determine a
coboundary relation between the non-abelian cocycle terms $\nu_{ijkl}$
and $\nu'_{ijkl}$. A certain amount of twisting takes
place,  however, and the   extra terms
$\{\widetilde{\zeta}_{ij},\,g_{klm}\}^{-1}$ and 
$\{\widetilde{m}'_{ijk},\,\te_{kl}\}$ need to  be inserted
 in their proper
locations,
just as the factor $\{\widetilde{m}_{ijk},g_{klm}\}^{-1}$ was
necessary in order to formulate equation \eqref{3coceq}.
 Once more, an equation such as  \eqref{cobouneq2} cannot be
viewed in isolation from its companion equation \eqref{cobouneq1}.
 In addition, any arrow in either of the monoidal  categories
$\mathrm{Ar}(\gc)$ or $\gc$ must be replaced by the corresponding 
one which is sourced at the
identity, without changing its name.
 The following definition summarizes the previous discussion.
\begin{definition}
Let $(\lambda_{ij},\, \widetilde{m}_{ijk},\,g_{ijk},\, \nu_{ijkl}) $
 and $ (\lambda'_{ij},\, \widetilde{m}'_{ijk},\, g'_{ijk},\, \nu'_{ijkl})$
 be a pair of 1-cocycle quadruples
 with values in the  weak crossed
square
\eqref{crsq}. These two cocycle quadruples are cohomologous if there
exists   a
quadruple $(r_i,\, \widetilde{\zeta}_{ij},\,\te_{ij},\,  b_{ijk})$
with values in  the weak crossed cube \eqref{crsq}. More precisely,
 these elements take their values in the square
\begin{equation}
\label{crsq3}
\xymatrix{
(\mathrm{Ar}_I\,\gc)_{U_{ijk}} \ar[d]_t \ar[r]^(.4)j 
& (\mathrm{Ar}_I\,\mathcal{E}q(\gc))_{U_{ij}}
\ar[d]^{t}\\
(\mathrm{Ob} \, \gc)_{U_{ij}} \ar[r]_(.4)j & 
(\mathrm{Ob}\, \mathcal{E}q(\gc))_{U_{i}} 
}
\end{equation}
in the positions 
\begin{equation}
\label{crsq4}
\left(
\begin{matrix}
 b_{ijk} & \widetilde{\zeta}_{ij}\\
\te_{ij} & r_{i}
\end{matrix}
\right) 
\end{equation}
The arrows $b_{ijk}$ and $\widetilde{\zeta}_{ij}$ are 
respectively of the form
\[\xymatrix@C=40pt{ I \ar[r]^(.23){b_{ijk}} & 
\lap_{ij} (\te_{jk}) \, \te_{ij}\, r_i(g_{ijk})\,
\te_{ik}^{-1}\,
(g'_{ijk})^{-1}}
\]
and
 \[
\xymatrix@C=40pt{
I  \ar[r]^(.3){\widetilde{\zeta}_{ij}} &  j(\te_{ij})\, r_i\, \lam_{ij}\,
  r_j^{-1}
\,{\lap_{ij}}^{-1}
}\]
and are required to 
satisfy the equations \eqref{cobouneq1} and \eqref{cobouneq2}. The
set of equivalence classes of 1-cocycle quadruples \eqref{crsq2},
 for the equivalence  defined by
these coboundary relations, will be 
 called  the \v{C}ech degree 1 cohomology
set for the open covering  $\mathcal{U}$ of $X$
with values in the weak crossed  square \eqref{crsq}. Passing to
 the limit over the
families of such open coverings of $X$, one obtains  the \v{C}ech
degree  1 cohomology
set of $X$ with values in this square.
 \end{definition}

\begin{remark}
{\rm If we consider as in Remark 6.6 a pair of normalized cocycle
  quadruples, 
the corresponding normalization  conditions on the coboundary terms  
\eqref{crsq4} are
\[
\begin{cases}
\te_{ij} = 1 \ \text{and} \ \,\widetilde{\zeta}_{ij} = 1 \
\text{when}\ i=j\\
b_{ijk}=1 \ \text{whenever} \ i=j\ \text{or} \ j=k
\end{cases}
\]}
\end{remark}
\noindent The discussion in paragraphs 
 6.2-6.4 can now be entirely summarized as
follows:
\begin{proposition}
The previous constructions associate to a $\gc$-2-gerbe $\pc$ on a
space $X$ an element of the \v{C}ech degree 1 cohomology
set of $X$ with values in the square \eqref{crsq}, and this element is
independent of the choice of local objects and arrows in $\pc$.
\end{proposition}
\noindent We refer to chapter 5 of \cite{2-gerbe} for the converse
 to this proposition, which
asserts that to each such 1-cohomology class corresponds a
$\gc$-2-gerbe, uniquely defined up to equivalence.

\begin{remark}
{\rm As we observed in footnote \ref{paracomp},
the proposition is  only true as stated  when the space $X$
satisfies an additional assumption such as paracompactness.
 The general case is discussed in \cite{2-gerbe},
where
 the  open
covering $\mathcal{U}$ of $X$ is  replaced by a hypercover.
 }
\end{remark}

 \end{document}